\documentclass[10pt,onecolumn,twoside,notitlepage,aps,pra,showpacs,longbibliography,nofootinbib]{revtex4-1}

\usepackage[hmargin=2.5cm]{geometry}
\usepackage[colorlinks=true,citecolor=blue,urlcolor=blue]{hyperref}
\usepackage{float}
\usepackage{graphicx}
\usepackage[sc,osf]{mathpazo}
\usepackage{amssymb,amsthm,amsfonts,amstext,amsmath}
\usepackage{color}
\usepackage{url}
\def\be{\begin{equation}}
\def\ee{\end{equation}}
\def\bea{\begin{eqnarray}}
\def\eea{\end{eqnarray}}
\def\bma{\begin{mathletters}}
\def\ema{\end{mathletters}}

\def\q0{\underline{0}}

\def\one{\leavevmode\hbox{\small1\normalsize\kern-.33em1}}

\newtheorem{theo}{Theorem}

\begin{document}

\title{Theoretical research without projects}
\author{Miguel Navascu\'{e}s and Costantino Budroni}

\affiliation{Institute for Quantum Optics and Quantum Information (IQOQI) Vienna\\ Austrian Academy of Sciences}

\begin{abstract}
We propose a funding scheme for theoretical research that does not rely on project proposals, but on recent past scientific productivity. Given a quantitative figure of merit on the latter and the total research budget, we introduce a number of policies to decide the allocation of funds in each grant call. Under some assumptions on scientific productivity, some of such policies are shown to converge, in the limit of many grant calls, to a funding configuration that is close to the maximum total productivity of the whole scientific community. We present numerical simulations showing evidence that these schemes would also perform well in the presence of statistical noise in the scientific productivity and/or its evaluation. Finally, we prove that one of our policies cannot be cheated by individual research units. Our work must be understood as a first step towards a mathematical theory of the research activity.\end{abstract}

\maketitle

\section{Introduction}

The introduction of performance-based funding at the end of the twentieth century, both at the level of institutions and single researchers, has had a significant impact on the current organization of the academic world. At the institutional level, research evaluation practices with the goal of distributing public funding have been introduced in several countries~\cite{Hicks2012}. At the same time, universities and research institutions rely more and more on grants to pay salaries, which are obtained by single researchers or consortia via a competitive procedure, in parallel to an increase of the percentage of non-permanent positions and competition among scientists \cite{Huisman2002academic, afonso2016, Kwiek2015changing,  AfonsoBlog}. The goal of these policies is to improve the performance of the research system, e.g., at a national or European level. On the one hand, these policies are supposed to move funds from the ``less efficient'' to the ``more efficient'' institutions and researchers, thus optimizing the fund allocation; on the other hand, the competition among researchers is supposed to drive their productivity.

At present, most research policies rely on the unfortunate idea that there is just one way of conducting science. In fact, there are two: experiments and theory, and science cannot advance without either. Theoretical researchers propose mathematical models to predict the outcomes of future experiments; and their experimental counterparts test the validity of such models in the lab. The programs currently used by grant agencies to promote research seem to be solely conceived with experimental science in mind. Not surprisingly, they are quite unfit for grant applicants from the theoretical sciences, such as mathematicians, computer scientists and theoretical physicists.

Let us explain why. In order to apply for a grant, most funding agencies demand a research project, i.e., they require scientists to detail their research activities within a period of $2$ to $5$ years in the future. While this scheme reflects the way many experimentalists organize their research agenda, it is incompatible with the current practice of theoretical research. One cannot ``plan'' discovering mathematical calculus, quantum cryptography or neural networks. Quite the opposite, some of the most celebrated ideas in history have arisen during the course of an unrelated investigation (see, e.g., \cite{Gillies2008book} for a few examples). It is by following these new threads that theorists keep their scientific productivity. On the contrary, stubbornly sticking to a single research line no matter what is a sound predictor of scientific sterility.

Agencies also demand theorists to lay out their ``methodology''. Namely, they expect theorists to explain how they intend to prove this or that theorem. The honest answer is that they don’t know. If they did, the theorem would be proven already, and they would not be applying for funds to crack it.

For the working theorist, applying for research funds is therefore a long and unethical task. It involves concocting an elaborate fantasy where the theorist pretends to know what theorems he/she will be proving in the next few years and through which particular mental processes. This activity, not clearly correlated with the applicant's success, takes a lot of time away from research \cite{Hippel2015_time_grant}. The product of these efforts, the project proposal, has no value whatsoever for society, and yet it is kept secret on the grounds of avoiding plagiarism. This lack of transparency makes grant panels unaccountable of any decision they make. 

Let us stress that the current grant system pushes theorists to lie in order to get funded. Indeed, if theorists carry their research sensibly, there will invariably be a mismatch between the original goals of the project and the final research output. So far no major scandal has transpired because the evaluators of a grant's final reports are researchers themselves: acknowledging that the system is flawed, they almost always award a positive assessment. This state of affairs, though, could change from one day to another, making thousands of theoretical researchers liable to a civil suit for fraud. In this direction, Gillies documents grant rules which advocate for the punishment of researchers who do not achieve the project goals \cite{Gillies2008book}. E.g.: forbidding them to apply again for the next two years, or reporting them to the head of their research institution.

We have reached this situation because, up to now, research policies have been based more on political fashion than on solid science. To progress beyond this point, we need an open scientific debate on research funding practices, where the scientific method is applied to the problem, i.e., with hypothesis, models, and experiments \cite{Azoulay2012}. The problem of research funding can be roughly divided into two sub-problems: first, the identification of the best evaluation method and corresponding parameters, e.g., of productivity or impact, and then the problem of maximization of such parameters given the available financial resources. 

Regarding the first problem, for the reasons above we believe that, at least for theoretical sciences, agencies and institutions should focus on funding people rather than projects \cite{Ioannidis2011comment}. That is, the evaluation of a researcher should be based on past scientific merits, as opposed to megalomaniac delusions.

The gross of the present paper addresses the second problem. Namely, presuming the existence of a quantitative measure of research productivity, estimated through the analysis of recent scientific activity, we investigate practical methods to optimize the total production of a global research system. 

We start by modeling the research system as a collection of agents or research units. A research unit can represent an individual scientist, a research group or a whole research institute. Each research unit possesses a ``scientific productivity function'' that relates how much science a given research unit can produce with the funds it holds to conduct research. We allow productivity functions to be probabilistic and time-dependent. They are also unknown, i.e., neither the research agency nor the scientists themselves can tell how they look like.

Relying on our mathematical model of the research activity, we show that there exist systematic procedures to decide the budget distribution at each grant call with the property that the total productivity of the research community will be frequently not far off its maximum value.

The simplest of such procedures is what we call \emph{the rule of three}, by which the funds $x_i^{k+1}$ for research unit $i$ after grant call $k+1$ are proportional to the research output $g_i^k$ of the unit during the $k^{th}$ term. If the total budget for science during the $(k+1)^{th}$ term is $X$ euros, this means that 

\be
x_i^{k+1}=X \frac{g_i^k}{\sum_jg_j^k}. 
\ee
\noindent The returns of this policy must be contrasted with those of ``excellence'' schemes, whereby, under equal research outcomes, researchers which were funded in the past have a greater chance of receiving further funds. Such policies can be shown to converge to configurations where the total scientific productivity is an arbitrarily small fraction of the maximum achievable by the research system. They are therefore better to be avoided.

We also study to what extent research policies can be cheated by dishonest research units. We conclude, for example, that hacks of the rule of three would require either influencing the evaluation stage or a coalition of research units. They are hence unlikely.

An important aspect of the current funding system is the fact that the increased competition and instability generate pressure among researchers, the so-called “publish or perish” culture, with possible negative consequences discussed in the literature, such as the focus on popular topics, short-term goals, and conservative research \cite{RzhetskyPNAS2015cons_reas,Wang2017bias,Gillies2008book} and the proliferation of useless research or even dishonest  practices \cite{Steen2013retract,Necker2014,RetractionWatch,Herteliu_2017}. 

The new funding framework that we advocate for is probably not a solution to the above problems, which are also closely connected to evaluation practices. Our framework, however, does not force theorists to engage in unethical practices, it is transparent and does not require the applicant to waste months of working time in writing project proposals. In addition, our mathematical analysis of scientific populations suggests that our grant schemes are relatively free from the so-called Matthew effect (i.e., ``the rich get richer and the poor get poorer.'')~\cite{Bol2018matthew}.

The paper is organized as follows. In Sect.~\ref{nProject}, we will introduce and motivate the use of a funding scheme not based on project proposals, but rather on the evaluation of recent-past productivity of single scientists or research institutions. In Sect.~\ref{prodFunc}, we will discuss which mathematical properties an idealized productivity function should possess. In Sect.~\ref{allocProb} we will define mathematically the problem of maximizing the total productivity of the system, given constraints on the total funding, and discuss possible ways of solving it, assuming the knowledge of all the parameters of the problem. In Sections \ref{fundPol}, \ref{moreCompl} we will adapt the previous analysis to the more realistic case of unknown model parameters, and explain how to extract a funding policy in this situation. In particular, we will perform numerical simulations to compare the performance of the different funding schemes under noise. In Sect. \ref{dishonest} we will analyze the security of one of our policies against dishonest players. In Sect.~\ref{discuss}, we will discuss possible extension of the model to handle more complex situations, e.g., competing funding agencies. Finally, we will present our conclusions.

\section{Funding policies not based on research projects}
\label{nProject}

The purpose of public research agencies is to help scientists generate useful knowledge, and the problem they face in each grant call is how to make sure that their funding reaches those in a position to generate such knowledge. 
Most funding agencies assign grants with a competitive process based on a peer-review evaluation of research proposals. As we argued at the introduction, this approach is not suitable for theoretical sciences because theorists cannot predict their future work activity.

A more appropriate indicator of the quality of future theoretical research is the recent past performance of the grant candidate. This motivates an alternative grant scheme for theoretical sciences based on the principle that, if candidates have recently shown a remarkable scientific productivity, it is worth funding them for the next few years so that they keep doing their good work. But what does ``good work'' mean?

There are many ways to quantify scientific productivity, and deciding which one suits best reflects a political stance. All such approaches fall in one of two main categories, namely, those based on bibliometrics and those based on peer-review. Bibliometric data has the advantage of being easy and cheap to obtain, e.g., through online databases containing publications and citations data, hence, it is often preferred by managers and administrators. However, there has been a proliferation of bibliometric indicators of scientific productivity and impact \cite{Waltman2016review,Abramo2017review}, often without a clear understanding of their pros and cons, from the perspective of evaluation and decision making. It has being argued that many of these indicators reflect ``what can be easily counted, rather than what really counts''\cite{Abramo2017review}. 

From the point of view of researchers, some methods may be considered fairer, such as expert assessment of the most important recent papers of the candidate. However, peer-review may be impractical in terms of cost and time, and  even be partially flawed (see, e.g., \cite{Birukou2011_alt_pr,Marsh2008improving} and references therein) showing, for instance, low reliability in the evaluations \cite{Jayasinghe2003multilevel} or low ratings for highly novel ideas \cite{creative}. Several authors have investigated whether peer-review and bibliometrics can be used together and how much they agree \cite{BacciniDeNicolao2016,Waltman2018_agree_pr,Campbell2010Bib_peer}. Bollen \emph{et al.} go further and propose an alternative funding scheme whereby the evaluation is conducted by the whole scientific community \cite{Bollen2014,Bollen2018}. 

In summary, notwithstanding the growing interest in the problem of research evaluation and the important results achieved so far, we feel that there is no general agreement on what the best evaluation methods are. As a consequence, we will leave this problem open and just start from the assumption of the existence of an abstract indicator of ``scientific productivity''.

Once a measurable figure of scientific productivity is established, the question is how to decide how much funds each research unit should receive. This is the problem we tackle in the rest of this paper. Curiously enough, we find that, given an agreed measure of scientific productivity, reaching an optimal allocation of research funds is not a political problem, but a mathematical one. In fact, we will show that under ideal conditions there exists a systematic procedure to decide the budget distribution at each grant call with the property that the total productivity of the research community will be frequently close to its optimal value.

At this point, it is important to remark which problems we are not addressing in our work and what possible use we can recommend or discourage. First, given the possible negative consequences of the publish-or-perish culture and the attitude towards experimenting with research policies, discussed in the previous section, we believe only a fraction of the entire research budget, e.g., at a national level, should be assigned through competitive grants. It is still an open question how much competition is desirable in academia, see, e.g., the discussion in Ref.~\cite{SANDSTROM2018} and references therein. Second, we leave open the question of at which level the evaluation and funding distribution should be applied. We will speak, generically, of ``research units''. Each research unit could be a single scientist, a research group, a department, a small research institution, or a university. We will occasionally speak of ``a scientist'' to provide examples and motivations for our assumptions, but the results of our work are independent of this choice. Third, we would like to remark that the methods and computational tools presented in this work are intended to aid human decisions. We do not advocate for a scenario where scientists are constantly evaluated by an algorithm that decides and directly modifies their salary. Finally, it is important to remark that political decisions are sometimes disguised as technical or scientific ones, e.g., budget cuts for universities and research institutions may be justified as technical decisions for the optimal use of available financial resources. The distinction between technical and political decisions should be made as clear as possible. We hope that separating the problem of funding from that of evaluation may bring clarity to the political debate.

\section{Scientific productivity functions: definition and properties}
\label{prodFunc}

\begin{figure}[b]
  \centering
  \includegraphics[width=12 cm]{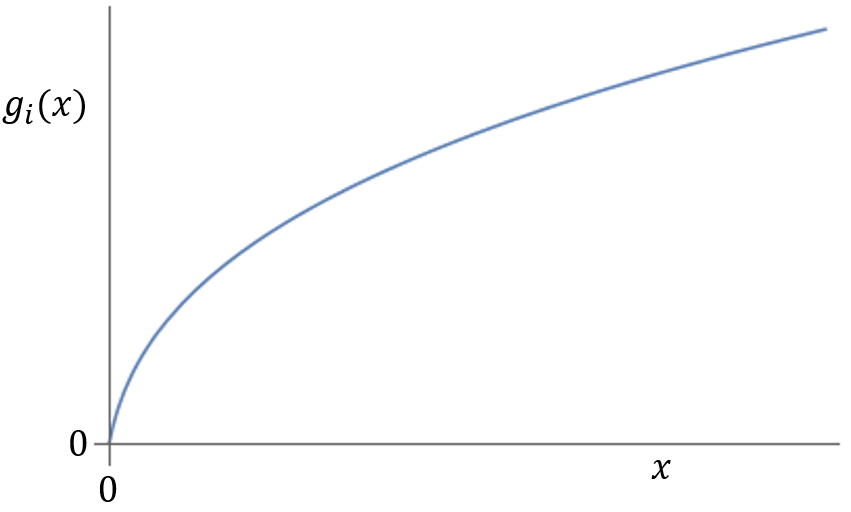}
  \caption{Expected productivity $g_i$ of a research unit as a function of its budget $x$. This picture illustrate the three main assumptions: the function is zero for zero budget,  it is non-decreasing, and it is concave.}
  \label{fig:productivity}
\end{figure}

Consider an idealized scenario where there is just one grant agency administering all public funds for research, and $N$ ``players'' (using a game-theoretic terminology) or research units apply for funding in consecutive grant calls. For further simplicity, we will adopt first a simple model where each player $i=1,...,N$ has a time-independent productivity function $g_i$. That is: if we award a player $x$ euros and demand it to use these funds within a time span $T$, then the scientific productivity of this player after time $T$, however we measure it, will be $g_i(x)$. Moreover, this quantity will be the same, independently of \emph{when} we awarded the player the research funds. In Sect.~\ref{moreCompl}, we will  present a generalization to probabilistic and time-dependent productivity functions.

\begin{figure}[b]
  \centering
  \includegraphics[width=15 cm]{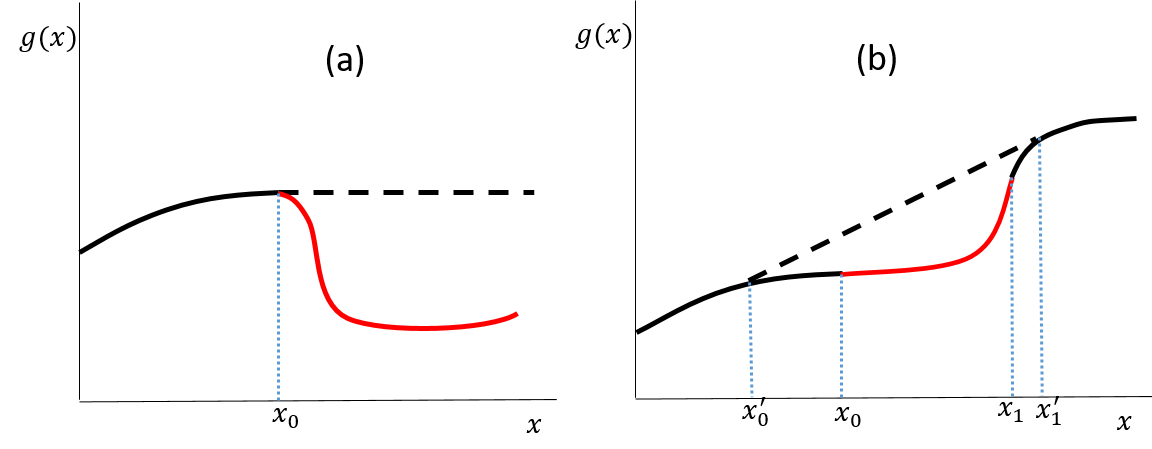}
  \caption{The productivity of a research unit as a function of the funding $x$. As a function of $x$, the scientific productivity $g$ of a rational player cannot have: (a) decreasing regions; or (b) convex regions. In both cases, using the same budget, the research unit can switch to a more favorable productivity function (dashed line) that is increasing and concave.}
  \label{prodEng}
\end{figure}

We have three main assumptions about the productivity functions $g_i$, which we will further discuss and motivate in the following. Namely,
\begin{itemize}
\item[(a)] $g_i(0)=0$, i.e., the productivity is zero when the budget is zero;
\item[(b)] $g_i$ is non-decreasing, i.e., if we increase the budget we should not decrease the productivity;
\item[(c)] $g_i$ is concave, i.e., the slope of the function is not increasing.
\end{itemize}
An example of a productivity function satisfying (a)-(c) is presented in Fig.~\ref{fig:productivity}.

Assumption (a) is quite straightforward: we do not expect a scientist without a salary to produce any science. There are indeed examples of outstanding individuals, such as Einstein, Bose or Gosset (Student), who carried out important theoretical contributions while working outside the academic world. However, all those individuals had also a salary, i.e., they had a monthly money input $x$ to play with. 

For assumption (b), we expect that if the player behaves rationally, the productivity function should not decrease with $x$. Indeed, suppose for example that $g(x_0)>g(x)$, for $x>x_0$, see Fig \ref{prodEng}. Then the player awarded with $x> x_0$ could simply spend $x-x_0$ euros to organize a conference and use the remaining $x_0$ euros to fund its research. Effectively, the player would then be operating according to a new increasing productivity function $\tilde{g}(x)$, with the property $\tilde{g}(x)\geq g(x)$ for all $x$.

Similarly, one can argue that any productivity function must be approximately concave, i.e., assumption (c). For suppose that, on the contrary, the increasing function $g(x)$ is convex in the region $[x_0,x_1]$, see Fig \ref{prodEng}. Fix $x_0'\leq x'_1$ and consider the following research strategy: if the funding $x$ satisfies $x\not\in [x'_0,x_1']$, then the scientist conducts research as usual, i.e., it will produce an output $g(x)$. If, on the contrary, $x\in [x'_0,x_1']$, then there exists a number $0\leq\lambda\leq 1$ such that $x=\lambda x'_0+(1-\lambda) x'_1$. In this case, we require the scientist to spend $x'_0\lambda$ euros for a fraction $\lambda$ of the total duration of the grant; and $x'_1(1-\lambda)$, for the remaining time $1-\lambda$. Assuming that, under a constant monthly salary, scientific productivity is time independent (namely, a scientist working for $2t$ months will produce twice as much as the same scientist working for $t$ months under the same salary), the total productivity will be $\lambda g(x'_0)+ (1-\lambda) g(x'_1)$. As shown in Fig \ref{prodEng}, one can choose $x_0'\leq x_0, x'_1\geq x_1$ such that the new effective productivity function $\tilde{g}(x)$ the scientist is operating under is concave. Moreover, $\tilde{g}(x)\geq g(x)$ for all $x\geq 0$.

The only problem with the above argument is that scientific productivity is just approximately linear with time. Indeed, one cannot expect $1000$ postdocs to advance significantly a new research line if they just have one day to do so (and we are overlooking the fact that very few would accept being employed for such a short time!). Hence, if $x'_0$ and $x'_1$ are very distant and $x\ll x'_1$, the previous scheme is not realistic.

In the following, though, we will assume for simplicity that the individual productivity is a concave function. This may not be very accurate to model the activity of a single scientist, but should be a good approximation to assess the productivity of a large group or a whole research institute. In sum, the shape of function $g_i(x)$ is expected to be approximately of the form depicted in Fig \ref{fig:productivity}. 

Note that, from the conditions of concavity and $g_i(0)=0$, it follows that 

\be
\frac{g_i(x)}{x}\geq g_i'(x), 
\label{concavFund}
\ee
\noindent for $x\geq 0$. Indeed, compute a first-order Taylor expansion of $g_i(0)$ on $x$. That gives us

\be
0=g_i(0)=g_i(x)-xg_i'(x)+\frac{g''_i(c)x^2}{2},
\ee
\noindent for some $c\in [0,x]$. Since $g_i(x)$ is concave, its second derivative is smaller than or equal to zero \cite{refConvex}. It follows that the right hand side of the above equation is upper bounded by $g_i(x)-xg_i'(x)$. Eq. (\ref{concavFund}) will be extensively used throughout the paper. If the inequality in (\ref{concavFund}) is strict for $x>0$, we will say that the function $g_i(x)$ is \emph{curved at the origin}. Intuitively, this means that, for any $a>0$, the productivity function $g_i(x)$ is not a straight segment from $x=0$ to $x=a$.

A family of productivity functions satisfying all these properties and rich enough to model interesting grant scenarios is the one given by ``power functions'' of the form $g(x)=Ax^\alpha$, where $A>0$, $\alpha\in (0,1)$. This family was already considered in \cite{bigScience}, where an attempt was made to estimate the average productivity function of a group leader. Moreover, the same power function, with an exponent smaller than (but close to) 1, was obtained by analyzing total citation counts (across 26 scientific disciplines) versus funds (Higher Education expenditure on Research \& Development expressed in Purchase Parity Power dollars) for OECD countries \cite{Cimini2014}.

\section{The problem of fund allocation}
\label{allocProb}

Let $X$ be the total funding that the agency can award in a given grant call. The goal of the funding agency is to identify the distribution of funds that maximizes the research output, given upper and lower bounds of the form $X^{-}_i\leq x_i\leq X^{+}_i$ on each player's budget $x_i$. The upper bounds stem from both the unwillingness of the individual to coordinate a large research group/institution and/or the desire of the funding agency of not concentrating a large amount of research funds in the hands of a few players. The lower bounds $\{X^{-}_i\}_i$ could correspond to negotiated minimum budgets for each research institute or public servant. Through the rest of the paper, the set of constraints 

\begin{align}
&\sum_jx_j=X,\nonumber\\
&X^{-}_i\leq x_i\leq X^{+}_i, \mbox{ for } i=1,...,N
\end{align}

\noindent will be denoted the \emph{funding conditions}. In case $X^-_i=0,X^+_i=\infty$, for $i=1,...,N$, i.e., in case the only restriction of the individual budgets is that they are non-negative, we will speak of \emph{free funding conditions}. In case $X_i^-=0$, we will speak of \emph{capped funding conditions}.

Ultimately, any funding agency wants to solve the optimization problem:

\begin{align}
g^\star(X)\equiv\mbox{maximize }&\sum_{i=1}^Ng_i(x_i),\nonumber\\
\mbox{such that } &\sum_{i=1}^Nx_i=X,\nonumber\\
&X^{-}_j\leq x_j\leq X^{+}_j,j=1,...,N.
\label{opt}
\end{align}

Since the funding conditions define a convex set and the objective function $\sum_{i=1}^Ng_i(x_i)$ is concave, any local maximum of $\sum_{i=1}^Ng_i(x_i)$ is also a global maximum. In other words: independently of our current budget configuration $\{x_i^0\}$, one can always identify the direction towards the optimal productivity by exploring how the objective function grows locally. It is also easy to prove that, as a function of the total funding $X$, $g^\star(X)$ is also concave.

For free funding conditions and fully homogeneous productivity functions, i.e., $g_i=g_1$ for $i=2,...,N$, the best strategy turns out to be distributing the funding equally among the researchers, in order to exploit the greater initial gradient of their productivity functions. That is, the solution of the above problem is $Ng_1\left(\frac{X}{N}\right)$. If $g_1(x)$ admits a first derivative at $x=0$, for $N\gg 1$, the latter quantity tends to $g_1'(0)X$.

Unfortunately, scientists can have very different productivity functions. Consider a scenario where each scientist $i$ has a power productivity function 

\be
g_i(x)=A_ix^{\alpha_i}. 
\label{powerF}
\ee

In Appendix \ref{solGeom} it is shown that the maximal productivity of this scientific population is given by

\be
g^\star(X)=\sum_i A_i(\alpha_iA_i\mu(X))^{\frac{\alpha_i}{1-\alpha_i}},
\label{soluPower2}
\ee
\noindent where $\mu(X)$ is computed by solving the equation

\be
\sum_i (\alpha_iA_i\mu)^{\frac{1}{1-\alpha_i}}=X.
\label{soluPower1}
\ee

\subsection*{Example}
It is at this point instructive to try to apply this simple model to ``real data''. Of course, this example is only illustrative of certain peculiar properties of the productivity functions and of the funding model. For instance, for simplicity, we will measure productivity simply by counting the number of papers, which is clearly a terrible quantifier, which \textbf{we do not endorse}. Moreover, contrary to the method present in Sect.~\ref{fundPol}, the current example uses the assumption of the specific productivity function of Eq.~\eqref{powerF}, which is in general not necessary. Finally, we do not claim that the numerical values obtained are particularly realistic, as they are extracted from only two data points and we provide no statistical analysis. Let us first go through the details of the example and, then, discuss at the end.

Since 1996, FWF Austrian Science Fund's START program provides the successful applicant with a funding amount between 0.8 and 1.2 million euros, to be spent in six years \cite{START}. The elegibility requirements demand that the doctoral degree of applicants be completed no less than two years and no longer than eight years before the deadline for submission of applications. It is thus not unreasonable to assume that most successful candidates did not have any prior funding, other than their own postdoc salary, before receiving the START grant. We randomly selected six START awardees, all of which work either in theoretical physics or mathematics, and estimated their scientific individual productivities by counting their number of peer-reviewed published papers in the six years prior to the year of the award and also in the six next years. For each candidate $i$ that provided us with two productivity points $g_i^1, g_i^2$ for each candidate. Complemented with the two funding inputs $x_i^1=435,780.00$ euros (the salary of a Senior Postdoc in Austria for six years) and $x_i^2=1,200,000.00$ euros (the maximum START funding), we had enough information to infer $A_i,\alpha_i$ for each researcher. 

Just for the matter of illustration, we have adopted the number of publications as a figure of merit. Since all the considered researchers received their START grants between 2007 and 2011, one would expect them to have been raised in the culture of ``publish or perish''. It is therefore sensible that most of them dedicated a substantial amount of effort to maximize their publication number.

The parameters of the so-computed productivity function for each researcher are displayed in the table below:

\begin{table}[H]
\begin{center}
\begin{tabular}{|c|c|c|}
\hline
$i$ & $A$ & $\alpha$\\
\hline
$1$ & $0.1447$ & $0.35212$\\
\hline
$2$ & $0.0213$ & $0.44621$\\
\hline
$3$ & $4.9761$ & $0.14333$\\
\hline
$4$ & $0.1574$ & $0.41840$\\
\hline
$5$ & $0.0312$ & $0.43619$\\
\hline
$6$ & $0.0076$ & $0.62321$\\
\hline
\end{tabular}
\end{center}
\caption{Power productivity functions of a population of $N=6$ scientists.}
\label{population}
\end{table}

\noindent Note that all exponents $\alpha$ are between $0$ and $1$, in agreement with our assumption that productivity functions are increasing and concave.

\begin{figure}
\includegraphics[width=10 cm]{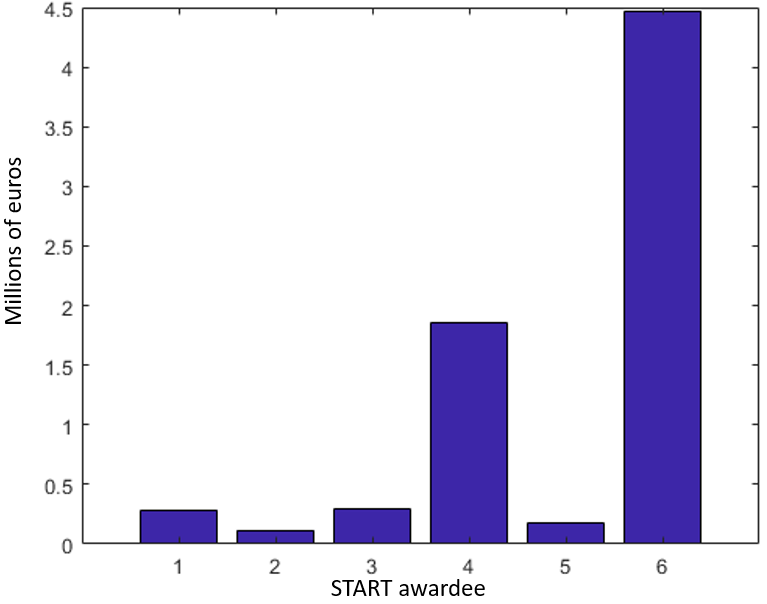}
  \caption{Optimal distribution of START funds over the $N=6$ researchers.}
  \label{xOptimal}
\end{figure}

The total amount of funds destined to these six researchers was $X=6\times 1.2$ million euros. Since the FWF distributed the funds equally among each researcher, the total productivity of this population due to the START grant is $\sum_{i=1}^6 g_i(X/6)\approx 184$. Using eq. (\ref{soluPower2}), however, we obtain a maximal productivity of $g^\star(X)\approx 225$. This is obtained by distributing budget $X$ in the way shown in Fig \ref{xOptimal}.

If the aim of the START program is to maximize the number of publications, then the program is just operating at $100\times 184/225\approx 81\%$ of its optimal yield. This percentage decreases as we increase the total funding $X$. Indeed, take $X=100$ million euros. In that case, an equal redistribution of the budget would produce an research output of $592$ publications. In contrast, the optimal allocation of such funds would give rise to $g^\star(X)\approx 886$. The performance of the egalitarian fund allocation would thus be $100\times 592/886\approx 67\%$. Of course, the aim of this example is to provide a simple illustration of our method, rather than to criticise the egalitarian fund distribution.

What is fundamental to notice is that, despite the extreme and unrealistic simplification of the model in this example, in particular the evaluation of ``productivity'' as a one-dimensional parameter, the player that obtains most of the funding in Fig.~\ref{xOptimal} is not necessarily ``the best researcher''. In fact, by changing the total amount $X$ available, one would change the optimal distribution of funds, hence, the (purported) induced ``ranking'' among researchers. More specifically, researchers with $\alpha_i\ll 1$ would receive most of the funds for low values of $X$, e.g., because they excel at working alone or in small groups, while researchers with $\alpha_i\approx 1$ would claim the greatest portion of the science budget in the high $X$ regime, e.g., because they excel at directing large research groups. This goes against the (possibly commonly accepted idea) that evaluation committees should choose, among a group of candidates for a grant, the one with the highest scientific productivity, irrespectively of the available resources. This concludes our example.
\\

In alternative to the estimation of the function parameters, one can apply a very intuitive rule of thumb to decide fund allocation. It consists of transferring funds between the different players until their average productivity rates $f_i(x_i)\equiv\frac{g_i(x_i)}{x_i}$ are as close as possible. That way, we arrive at a final distribution of funds $\{x^\sharp_i\}_i$ such that, for any $i\not=j$, $i,j=1,...,N$,

\be
f_i(x^\sharp_i)> f_j(x^\sharp_j) \mbox{ implies that } x^\sharp_j=X^{-}_j \mbox{ or } x^\sharp_i=X^{+}_i.
\label{condZero}
\ee

In principle, one can achieve such a configuration by solving the optimization problem:

\begin{align}
\mbox{maximize }&\sum_{i=1}^N G_i(x_i),\nonumber\\
\mbox{such that } &\sum_{i=1}^Nx_i=X,\nonumber\\
&X^{-}_j\leq x_j\leq X^{+}_j,j=1,...,N,
\label{optZero}
\end{align}
\noindent where $G_i(x)\equiv \int_0^x dy \frac{g_i(y)}{y}$. It can be easily proven that, for each $i$, $G_i(x)$ is a concave function. This means that the maximization (\ref{optZero}), like (\ref{opt}), does not risk getting stuck in local maxima.

Define  $g^{\sharp}(X)=\sum_ig_i(x^\sharp_i)$, where $\{x^\sharp_i\}_i$ is the configuration of funds maximizing (\ref{optZero}). Since $\{x^\sharp_i\}_i\not= \{x^\star_i\}_i$, we expect that $g^\sharp(X)<g^\star(X)$, i.e., this manner of allocating funds will not be optimal in general. 

In this regard, consider a bipartite ($N=2$) scientific population where, for $x\in (0,1]$ player $1$ has an almost constant productivity function $g_1(x)= 1$ for $x>0$, while player $2$'s productivity function is linear, i.e., $g_2(x)= x$. Then the optimal funding configuration consists in assigning player $1$ an infinitesimal amount of funds ($x^\star_1=\delta$), while giving player $2$ the rest ($x^\star_2=X-\delta$). The supremal (not maximal) productivity is thus $g^\star(X)=1+X$. On the contrary, using the above rule of thumb, it is easy to see that, for $X<1$, $x_1^\sharp=X$, $x_2^\sharp=0$, and so $g^\sharp(X)=1$. As $X$ gets close to $1$, the fraction $g^\sharp(X)/g^\star(X)$ tends to $\frac{1}{2}$. Not only we do not achieve the maximal productivity, but the funding distributions in one case and the other are complete different!

Nonetheless, it is possible that, while not being optimal, $g^\sharp(X)$ is not that far off the optimal scientific productivity. In this line, we have the following result.

\begin{theo}
\label{egregius}
Consider a scientific population characterized by productivity functions $\{g_i(x)\}_{i=1}^N$, subject to capped funding conditions $\{0\leq x_j\leq X_i^+,\sum_ix_i=X\}$. Then

\be
g^\sharp(X)\geq \frac{1}{2}g^\star(X).
\label{super}
\ee
\end{theo}

\noindent In other words: even though the grant scheme (\ref{zero_order}) is suboptimal, its performance is, at worst, half of the optimal one. Moreover, due to the example above, the constant $\frac{1}{2}$ cannot be improved.
See Appendix \ref{proofApp} for a proof.

As shown in Appendix~\ref{app:cosmetic}, for free-funding conditions, and provided that the slope of each function $g_i$ at $x=0$ is ``big enough'' (more concretely: $f_i(0)>f_j\left(\frac{X}{N-1}\right)$ for all $i,j$), we have that the optimal configuration $\{x_i^\sharp\}_i$ will satisfy

\be
\frac{g_i(x_i^\sharp)}{x_i^\sharp}=\frac{1}{\lambda}, \text{ for all } i,
\label{limit}
\ee
\noindent for some $\lambda>0$. 

Note that the slope condition holds for power productivity functions (since $f_i(0)=\infty$ for all $i$). For populations of research units described by such functions, we can thus use eq. (\ref{limit}) to derive an explicit expression for $g^\sharp(X)$, namely:

\be
g^\sharp(X)=\frac{X}{\lambda(X)},
\label{zeroPower2}
\ee
\noindent where $\lambda(X)$ is obtained by solving the equation:

\be
\sum_i (A_i\lambda)^{\frac{1}{1-\alpha_i}}=X.
\label{zeroPower1}
\ee

Applying formula (\ref{zeroPower2}) to the scientific population described in Table \ref{population} of the example, with $X=6\times 1.2$ million euros, we obtain $g^\sharp(X)\approx 213$ publications. This represents an efficiency of $100\times g^\sharp(X)/g^\star(X)\approx 95\%$. For $X=100$ million euros, we obtain $g\sharp(X)\approx 861$, with an efficiency of $100\times g^\sharp(X)/g^\star(X)\approx 97\%$. This allocation scheme thus seems to give a good performance when applied to real scenarios.

\section{Funding policies}
\label{fundPol}

Unfortunately, neither the funding agency nor the scientists themselves know the explicit form of their productivity functions. So how can a funding agency expect to solve problem (\ref{opt})?

In the following two sections, we provide a number of \emph{automated} methods to carry out this task. Under some assumptions on scientific production, some of these methods are guaranteed, on their own, to steer the productivity of a scientific population near its maximum possible value. 

Nonetheless, \textbf{these computational tools are intended to be used by human agents as an aid to reach a final budget decision}. Note that the sole purpose of these tools is to maximize a given figure of merit, irrespective of any other considerations. We strongly doubt that the whole scientific enterprise can be reduced to an optimization problem. Thus, by removing human intervention completely from scientific policy decision-making, we risk reaching a dystopian scenario where any aspect of science other than an agreed objective function is viewed as an obstacle towards the maximization of the latter. On the other hand, we know, from the world of chess, that in some situations the best decision-makers are neither human nor artificial, but a team of both kinds of entities. We are therefore confident that funding agencies will greatly benefit from the ideas that we present next.

We will start by dividing time in \emph{terms}, of $s$ years each. At the end of the $k^{th}$ term, the funding agency announces the $(k+1)^{th}$ \emph{call}, and, after a proper evaluation, distributes the funds in such a way that player $i$ receives $x_i^{k+1}$ euros to be spent on the $(k+1)^{th}$ term. 

Now, let us forget for the time being that we ignore the objective function $g(\bar{x})\equiv\sum_{i=1}^Ng_i(x_i)$. A very effective tool to solve maximization problems like (\ref{opt}) is the projected gradient method \cite{subgradient}. The output of this method is a sequence of feasible budgets $\bar{x}^1,\bar{x}^2,\bar{x}^3,...$ with the property that 

\be
\frac{1}{k}\sum_{j=1}^kg(\bar{x}^j)\approx g(\bar{x}^\star),
\label{orbits}
\ee
\noindent with $\bar{x}^\star$ being the optimal configuration. Each budget $\bar{x}^{k+1}$ is obtained from the previous one $\bar{x}^k$ by the following iterative equation:
\be
\bar{x}^{k+1}=\pi_B \left(\bar{x}^k+\epsilon\bar{\nabla}g(\bar{x}^k)\right),
\label{fundamentalGradMet}
\ee
where $\epsilon>0$ is a free parameter known as the \emph{learning rate} and $\pi_B(\bar{z})$ denotes the closest vector $\bar{y}$ to $\bar{z}$ (in Euclidean norm) belonging to the set $B=\{\bar{x}:\bar{x}\geq 0$, $\sum_{i=1}^Nx_i =X, 0\leq x_j\leq X_j,\forall j\}$ of allowed budget configurations. Computing $\pi_B(\bar{z})=\mbox{arg}\min_{\bar{y}\in B}\|\bar{z}-\bar{y}\|_2$ can be cast as a semidefinite program \cite{sdp}, a type of optimization problems which we know how to solve efficiently.

Now suppose that, at the $(k+1)^{th}$ call of the grant, we chose to distribute the funds according to $\bar{x}^{k+1}$ in eq. (\ref{fundamentalGradMet}), and that we repeated this operation in all subsequent calls. Since $g(\bar{x})$ is (approx.) concave, by eq. (\ref{orbits}), very frequently the current budget distribution would (approx.) maximize the total scientific output of the community.

Our problem is, however, that we don't know $\bar{\nabla} g(x)$. Consider then the following modification of the iterative equation:

\be
\bar{x}^{k+1}=\pi_B \left(\bar{x}^k+\epsilon \tilde{\nabla}^kg\right),
\label{new_algo}
\ee

\noindent where we have approximated the gradient 

\be
\bar{\nabla} g(\bar{x}^k)=\left(\frac{\partial g_1(x^k_1)}{\partial x^k_1},\frac{\partial g_2(x^k_2)}{\partial x^k_2},...,\frac{\partial g_N(x^k_N)}{\partial x^k_N}\right)
\ee

\noindent by the vector 

\be
\tilde{\nabla}^kg=\left(\frac{g_1^k-g_1^{k-1}}{x_1^k-x_1^{k-1}},\frac{g_2^k-g_2^{k-1}}{x_2^k-x_2^{k-1}},...,\frac{g_N^k-g_N^{k-1}}{x_N^k-x_N^{k-1}}\right),
\ee
\noindent where $g_i^k$ is the declared scientific production of player $i$ at the end of the $k^{th}$ term (that should equal $g_i(x^k_i)$, if the player is being honest, see Section \ref{dishonest}).

One can show that the iteration scheme (\ref{new_algo}) also satisfies eq. (\ref{orbits}), see Appendix \ref{proofConv}. Note, though, that the vector component $(\tilde{\nabla}^kg)_i$ can be computed given the funding $x^{k-1}_i, x^{k}_i$ received by the candidate in the last two grant calls and the corresponding research outputs $g_i^{k-1}, g_i^{k}$. Hence this procedure can be implemented in practice. By using eq. (\ref{new_algo}) to decide the budget distribution in the $(k+1)^{th}$ term, we make sure that, in the long run, research funds are distributed in an optimal way. Note that there may be situations where we lack data to compute $\tilde{\nabla}^kg(x)$, e.g.: the candidate just finished the PhD studies, or had a child-raising break. In those cases, one can replace $(\tilde{\nabla}^kg)_i$ by $\frac{\hat{g}_i}{\hat{x}_i}$, where $\hat{g}_i, \hat{x}_i$ are, respectively, the last known scientific production of the candidate and the science funds it was enjoying at the time.

The recursive method (\ref{new_algo}), that we will in the following refer to as the \emph{gradient scheme}, is an instance of a \emph{grant policy}. There are others. Consider, for instance, the following one:

\be
x_i^{k+1}=\pi_B\left(x_i^k+\epsilon \frac{g_i^k}{x_i^k}\right).
\label{zero_orderGrad}
\ee

This is none other than the gradient method, applied to solve optimization problem (\ref{optZero}). For $\epsilon$ small enough, the orbit $(\bar{x}^k)_k$ will often be very close to the optimal configuration $\{x_i^\sharp\}_i$. Moreover, for capped funding conditions, Theorem \ref{egregius} guarantees that the corresponding total productivity will be, at least, one half of the optimal one. Note that this scheme only requires knowledge of the total scientific budget $X$ and the immediate past performance of the players: it is a \emph{zero-order scheme}, as opposed to the first-order scheme (\ref{new_algo}), that requires information of the last two grant calls. We will dub this policy the \emph{average rates scheme A}.

Alternatively, one can use (for free funding conditions) the iterative method:

\be
x_i^{k+1}=X\frac{g_i^k}{\sum_{j=1}^N g_j^k}.
\label{zero_order}
\ee
\noindent 
This is also a zero-order scheme, where we do not even need to know the funds $\{x^k_i\}_i$ awarded to each researcher at call $k$ in order to decide the funding distribution of call $k+1$.

Interestingly, if the initial distribution of funds $\{x_i^0\}$ satisfies $x_i^0>0$ for all $i$ and all productivity functions are curved at the origin, it can be proven that this policy converges exponentially fast to the configuration $\{x_i^\sharp\}_i$ of Eq.~\eqref{limit} (see Appendix~\ref{app:conv_zero} for a proof). We will refer to this grant policy as the \emph{average rates scheme B}, or, more colloquially, as \emph{the rule of three}, since, given $x^{k+1}_1, g_1^k$, any other player $i=2,...,n$ can use $g_i^k$ and the rule of three to compute its future funding.

Finally, there is another grant policy, that we will hereby call the \emph{standard scheme}, by which the funding of each researcher is proportional to the funds it received in the previous grant call and its productivity. That is:

\be
x_i^{k+1}=\pi_B\left(X\frac{x_i^kg_i^k}{\sum_j x_j^kg_j^k}\right),
\ee
\noindent where $\pi_B$ denotes, as before, the projection onto the set of budgets satisfying the funding conditions. The standard scheme reflects the growing perception in the theoretical physics community that the probability of being awarded a grant for a theoretical project grows with both the productivity of the candidate and the funding obtained in the past. For an evidence, consider the following extract from the Application Guidelines for Stand-Alone Projects, FWF Austria Science Fund: \emph{``Most important research projects funded in the past (no more than 5). [...] For each project, please provide the following information: Project title, funding agency, project duration (from/to) and amount of funding granted''}.

Both the distribution of funds and the final total productivity of this scheme depend significantly on the initial budget configuration $\bar{x}^0$. Indeed, consider free funding conditions, and assume that the productivity functions are all identical and equal to the power function (\ref{powerF}). Then one can show that the standard scheme converges to a fund configuration where just the players $i$ with maximum values of $x^0_i$ receive any funds whatsoever. For generic initial configurations $\bar{x}^0$, only one player $i$ will satisfy this demand, in which case the asymptotic total productivity of the standard scheme will be $g^\flat(X,\bar{x}^0)=AX^\alpha$. This has to be compared to the optimal productivity $g^\star(X)=N^{1-\alpha}A X^\alpha$, achieved by the ``egalitarian'' configuration $x_i=\frac{X}{N}$, for $i=1,...,N$. In this example, the final configuration enforced by the standard scheme is thus maximally unfair, with just one player holding all the resources, instead of an equal distribution of funds. In addition, for large populations of scientists ($N\gg 1$), the quotient $\frac{g^\flat(X,\bar{x}^0)}{g^\star(X)}$ becomes arbitrarily small. This gives some theoretical grounds for Nobel Laureate Jeffrey C. Hall's remarks in \cite{cell}: 

\begin{quotation}
\emph{I can't help feel that some of these [scientific] `stars' have not really earned their status. I wonder whether certain such anointees are `famous because they're famous'. So what? Here's what: they receive massive amounts of support for their research, absorbing funds that might be better used by others.}
\end{quotation}

\section{Probabilistic time-dependent productivity functions}
\label{moreCompl}

In realistic scenarios, it is expected that a player's productivity will not only depend on funding, but also on a number of variables which escape our control (health, lack of sleep, love affairs...). We can model the effect of these variables by postulating that productivity functions must be probabilistic. Actually, in the real world things are even more complicated: productivity functions vary with time, as researchers acquire new knowledge and skills, or their motivations waver. In these conditions formula (\ref{new_algo}) is not guaranteed to generate orbits close to an optimal productivity. 

Suppose then that $g_i$ is a probabilistic function that varies with time, i.e., the productivity of player $i$ at the end of term $k$ is a random variable $g^k_i$ of the form $g_i^k=g_i(x,k)$. Then deciding which quantity to optimize in this scenario is again a political (subjective) matter. A reasonable figure of merit, that we will use from now on, is the average scientific productivity at each term $k$. Note that one can repeat the arguments in Section \ref{prodFunc} to suggest that $\langle g_i(x,k)\rangle$ should also be increasing and approximately concave in $x$. Our goal is therefore to identify a policy to decide the fund allocation $\bar{x}^j$ at each grant call $j$, such that, for $k\gg 1$, 

\be
\frac{1}{k}\sum_{j=1}^k g(\bar{x}^j,j)\approx \frac{1}{k}\sum_{j=1}^k \max_{\bar{x}\in B}\langle g(\bar{x},j)\rangle
\label{orbitsProb}
\ee
\noindent with high probability.

This puts us in a conundrum. On one hand, a single estimate of $g_i(x_i,t)$ does not allow us to assess its average value, which we need to know in order to maximize the average productivity of the whole community. On the other hand, we cannot rely on the early past history of the candidate, because the productivity function also changes with time.

One possibility is to apply the gradient scheme (\ref{new_algo}), but with a correction that guarantees that random fluctuations do not squander the optimum budget. Note that, using the grant policy (\ref{new_algo}), it could be the case that a candidate receives the same funding twice consecutively, $x^k_i=x^{k-1}_i$, but outputs different results $g_i(x^k_i)\not=g_i(x^{k-1}_i)$. That would lead us to estimate an infinite gradient that would either put all the future budget in the hands of this candidate, or reduce its budget to $0$ in the present grant call. In such a predicament, it is more convenient to use the corrected formula

\begin{align}
&z_i^{k+1}=x_i^k+\epsilon H(x_i^k,x_i^{k-1},g_i^k,g_i^{k-1}),\nonumber\\
& \bar{x}^{k+1}=\pi_B \left(\bar{z}^{k+1}\right),
\label{cosica}
\end{align}
\noindent where $H(x_i^k,x_i^{k-1},g_i^k,g_i^{k-1})$ is a ``filtered version'' of $(\tilde{\nabla} g)_i$, defined as:

\begin{align}
H(x_i^k,x_i^{k-1},g_i^k,g_i^{k-1})=&0, \mbox{ if } (\tilde{\nabla}^k g)_i<0,\nonumber\\
&\frac{g_i^k}{x^k_i},\mbox{ if }(\tilde{\nabla}^k g)_i>\frac{g^k_i}{x^k_i},\nonumber\\
&(\tilde{\nabla}^k g)_i, \mbox{ otherwise}.
\label{hache}
\end{align}
The filter's goal is to get rid of non-sensical estimations $\tilde{\nabla}^k g$ of the actual gradient $\bar{\nabla}^kg$ of the objective function. Indeed, since $g_i(x)$ is increasing, it can't be that $g_i'(x)<0$. Similarly, by eq. (\ref{concavFund}), $g'_i(x)\leq f_i(x)$.

For $g$ deterministic and time-independent, $H(x^k_i,x^{k+1}_i,g^k_i,g^{k+1}_i)=\tilde{\nabla}^kg$. The policy (\ref{cosica}) is, in this scenario, equivalent to algorithm (\ref{new_algo}), and so its outputs $(\bar{x}^{k})_k$ will satisfy eq. (\ref{orbits}). We leave as an open question under which conditions eq. (\ref{orbitsProb}) is satisfied in the probabilistic, time-dependent case. For the rest of the article, the use of formula (\ref{cosica}) to decide the funding allocation will be dubbed \emph{the modified gradient scheme}.

Alternatively, we can resort to the average rates schemes (\ref{zero_orderGrad}), (\ref{zero_order}). Since these policies only take into consideration productivity data and funds from the previous grant call, one would expect them to be even more robust against the time evolution of the productivity function. Moreover, by inspection of eqs. (\ref{zero_orderGrad}), (\ref{zero_order}) it is clear that small random fluctuations on the productivity will hardly affect the resulting budget configuration after a few calls. In fact, for the average rates scheme A, it can be proven (see Appendix \ref{convA}) that, if the productivity functions change slowly with time, then with high probability we have

\be
\frac{1}{k}\sum_{j=1}^k g(\bar{x}^j,j)\approx \frac{1}{k}\sum_{j=1}^k \langle g(y^j,j)\rangle.
\label{orbitsProbA}
\ee
\noindent Here $y^j$ denotes $\bar{x}^\sharp$, as defined in (\ref{condZero}), for the functions $g_i(x)\equiv \langle g_i(x,j)\rangle$. By Theorem \ref{egregius}, this means that, under capped funding conditions, the average rates Scheme A is guaranteed to produce on average at least half of the optimal productivity.

So far we have discussed three different grant policies. Which one shall we use in practice? To help us answer that question, we will next compare their performance in a number of numerical simulations. 

In each simulation, we will consider a population of researchers with time-independent, non-deterministic productivity functions. Their average productivity functions will be given by Table \ref{population}. To model both the random fluctuations in productivity and the volatility of scientific evaluation, we assign to each scientist a measure of unpredictability $0\leq U\leq 1$: the actual productivity of the scientist is taken to be $g_i(x)(1+u_i)$, where $u_i$ is a random number chosen uniformly from the interval $[-U,U]$. In our simulations, we studied three cases of interest: $U=0$ (no noise), $U=1/8$, (low noise) and $U=1/2$ (high noise). 

\begin{figure}
  \includegraphics[width=17 cm]{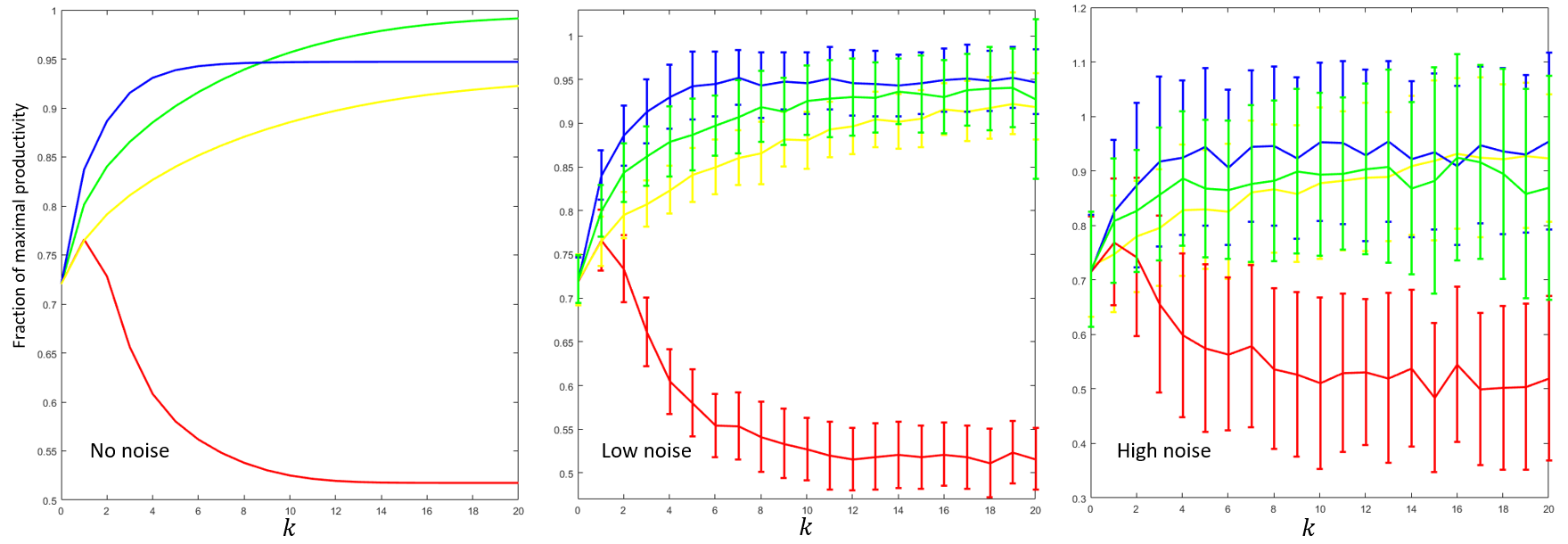}
  \caption{\textbf{Productivity as a function of the grant call $k$, for different scientific policies.} The colors green, yellow, blue and red denote, respectively, the modified gradient scheme, the average rates scheme A, the rule of three and the standard scheme. Starting from a random distribution of funds, we study the performance of different policies under increasing amounts of statistical noise. For both the modified gradient scheme and the average rates scheme A, we chose $\epsilon=\frac{0.25*X}{\max_i(g_i^0/x_i^0)}$. If each term lasts $s=4$ years, then, in all cases, the rule of three would require $12$ years to steer the community to a configuration where the average scientific production is greater than $90\%$ of the maximal value.}
  \label{policiesPic}	
\end{figure}

Starting with the random funding distribution $\bar{x}^0 =(1.6804, 1.8683,0.2619,1.8839,1.3043,0.2012)\times 10^6$, we estimated the fraction of the maximal productivity achieved via different policies at each call $k$. The normalized average productivity is depicted in Fig \ref{policiesPic}, together with its variance. As we can see, even under low statistical noise the rule of three performs slightly better than the modified gradient scheme after a reasonable number of calls ($5,6$), and substantially better than the standard scheme. In the asymptotic limit, the performances of the average rates schemes A and B are comparable, but the latter converges faster to the optimal value.

Which one of these policies should be chosen depends on the typical form productivity functions, a matter that can only be decided through experiment. If the relevant productivities exhibit almost no statistical noise, the modified gradient scheme shall be preferred. On the contrary, under free funding conditions and a fair amount of noise, the rule of three seems to be the wisest choice.

\section{Dishonest players}
\label{dishonest}

So far we have been assuming that all players are honest, namely, that they won't try to play the system to obtain more funds than they should. Such is a very naive position: ideally, we would like to have research policies which cannot be played. In this spirit, we will next study the security of the rule of three against dishonest participants.

For simplicity, we will carry the analysis under the assumption that the dishonest player, Daniel, belongs to a large scientific population where the budget distribution is almost stationary (i.e., it is close to an equilibrium). We will also assume that Daniel's optimal productivity function is time-independent, deterministic and curved at the origin. Finally, we will suppose that, if Daniel plays honestly, in the asymptotic limit his budget $x^\sharp_i$ will be greater than zero.

Following eq. (\ref{zero_order}), by producing an amount of science $g_i^k$ at the end of the $k^{th}$ term, Daniel will receive the funds $x_i^{k+1}=X\frac{g_i^k}{\sum_j g_j^k}$ at the $(k+1)^{th}$ call. Since the population Daniel belongs to is large and its funding distribution close to stationarity, $\frac{X}{\sum_j g_j^k}$ will hardly vary with $g_i^k$ and $k$, so we will take it constant, i.e., we will assume that $x_i^{k+1}=\lambda g_i^k$, for some $\lambda>0$. 

Suppose that Daniel plays honest. Then, given an initial amount of funds $x_i^0$, he will invest them all on research, thus producing $g_i^0=g_i(x_i^0)$ and earning $x_i^1=\lambda g_i(x_i^0)$ at the end of the process. Iterating, we find that the Daniel's funds will follow the orbit $x^0_i, \lambda g_i(x^0_i), \lambda g_i(\lambda g_i(x^0_i)), \lambda g_i(\lambda g_i(\lambda g_i(x^0_i))),...$. It can be shown that, in the limit, he will be receiving $x^\sharp_i>0$ after each call, with $x^\sharp_i$ defined by the relation $\lambda g_i(x_i^\sharp)=x_i^\sharp$. In other words, 

\be
\lim_{k\to\infty}\frac{1}{k+1}\sum_{j=0}^k x^j_i=x^\sharp_i.
\ee

Now, we let Daniel be dishonest. How could he get more than $x^\sharp_i$ on average? If the funding agency requires each player to spend all his/her funds by the end of the term, then the only thing Daniel can do is keep a fraction of his scientific results secret. Perhaps, by declaring a vast accumulation of scientific achievements in one go he may manage, on average, to squeeze more money out of the grant agency.

More specifically: starting with the funds $x_i^0$ and an amount of undeclared scientific results $g_i^{+}$, Daniel would produce an output $g_i(x_i^0)$. He could then declare, at the end of the term, that he has produced an amount of results $g_i^0$, with $0\leq g_i^0\leq g_i^++g_i(x_i^0)$. In turn, the grant agency would reward him with an amount $x_i^1=\lambda g_i^0$, that he would use to produce $g_i(x_i^1)$ results. In the next term, he would declare $g_i^1$, with $0 \leq g_i^1\leq g_i^++g_i(x_i^0)-g_i^0+g_i(x_i^1)$ results, and so on. Daniel's \emph{declaration strategy} is summarized in Table \ref{cheating}. 

\begin{table}[H]
\begin{center}
\begin{tabular}{|c|c|c|c|}
\hline
term & funds received & reported productivity &unreported productivity\\
\hline
$0$ & $x_i^0$ & $g_i^0$ & $g_i^++g_i(x_i^0)-g_i^0$\\
\hline
$1$ & $x_i^1$ & $g_i^1$ & $g_i^++g_i(x_i^0)-g_i^0+g_i(x_i^1)-g_i^1$\\
\hline
$2$ & $x_i^2$ & $g_i^2$ & $g_i^++g_i(x_i^0)-g_i^0+g_i(x_i^1)-g_i^1+g_i(x_i^2)-g_i^2$\\
\hline
\vdots & \vdots & \vdots & \vdots\\
\hline
\end{tabular}
\end{center}
\caption{Daniel's cheating strategy.}
\label{cheating}
\end{table}

The second column denotes the funds received at the beginning of the call. Columns $3$ and $4$ denote, respectively, the reported and unreported scientific production when the funds run out. Note that, for this table to represent a valid strategy, $x_i^{k+1}=\lambda g_i^k$ and the elements of column $4$ must be greater than or equal to $0$. The question is whether Daniel can choose $g_i^1,g_i^2,...$ in such a way that, on average, he will obtain more funds than acting honestly. That is, whether, for high $k$, $\frac{1}{k+1}\sum_{j=0}^{k} x^j_i>x^\sharp_i+\delta$, with $\delta>0$. In Appendix \ref{secApp} we show this not to be the case. 

Thus, in the long run, Daniel will not win anything by delaying the publication of his scientific production. On the contrary, by not publishing his results as soon as he produces them, he risks being scooped by some other player, in which case he would be losing well-deserved science funds. It is therefore in Daniel's interest to play honest and declare all his productivity by the end of each term.

Note that this security analysis required assumptions on both Daniel's productivity function (determinism, time-independence and curvature at the origin) and the overall behavior of the scientific population. It would be interesting to find out whether security also follows if such assumptions are dropped. Similarly, it would be interesting to see if the modified gradient scheme and the average rates scheme A are also secure under dishonest participants. 

Most crucially, we have not studied the possibility that the evaluation of the scientific production of the players is not impartial. Some studies suggest that for securying fund it is more important how researchers build their collaboration network than what publications they produce and whether they are cited  \cite{networking}. Actually, in some fields such a behavior has greatly influenced the distribution of research funds in the past \cite{trouble}. At the moment we do not have a solution for this problem, other than hoping that not so many scientists engage in this practice.

\section{Discussion}
\label{discuss}
Here we will examine some shortcomings of our model for research funding (\ref{opt}) and discuss how the latter can be improved for its use in real-world scenarios. This section is much more technical than the others and can be skipped on a first reading.

\subsection{Assumptions on the productivity function}
In section \ref{prodFunc}, we argued that the productivity function $g_i(x)$ of each player $i$ must be increasing and approximately concave. We did so by reasoning that any rational player who knows the shape of its productivity function can improve it piece by piece until it becomes increasing and approximately concave. The underlying assumptions are that the agent knows its productivity function, that it is interested in maximizing it and that it acts rationally. These three conditions may not be met in practice. 

Consider, for instance, the second one. Suppose that the goal of the funding agency were to maximize the total number of publications, while the personal goal of player $i$ is to maximize the quality of the said publications. Then, given more funding, player $i$ would not use it to increase its publication number, but to hire better researchers and hence produce better papers. In such circumstances, there may not be a simple relation between $x_i$ and the productivity $g_i$ measured by the agency. In principle, $g_i(x)$ could be decreasing, or convex.

The maximization of non-concave functions is a conventional problem in artificial intelligence, where the accuracy of the output of a neural network depends non-trivially on a number of continuous parameters. There exist a number of methods to achieve this effect, see \cite{compendium}. Unfortunately, all of them require a reliable estimate of the gradient $\bar{\nabla}^k g$. Under very low statistical fluctuations and productivity functions independent of time, we can approximate $\bar{\nabla}^kg$ by $\tilde{\nabla}^k g$ as in (\ref{new_algo}). In the general case, though, it is unclear what to do when one of the components of $\tilde{\nabla}^kg$ is negative or very high. Indeed, since $g_i$ may not be neither increasing nor concave, we cannot assume that neither $g'_i(x)\geq 0$ nor eq. (\ref{concavFund}) holds, and so we are not entitled to filter $\tilde{\nabla}^kg$ as in eq. (\ref{hache}).

Another tacit assumption in (\ref{opt}) is that the productivity $g_i$ of a player $i$ just depends on its funding $x_i$ (and not on the funding $\{x_j:j\not=i\}$ of all the other players). This condition does not capture frequent real-world situations where two or more research institutes compete for the same gifted group leader. A more realistic model for scientific productivity would posit that there exists a global productivity function $g(x_1,...,x_N)$ that does not necessarily decompose as a sum of independent productivities, i.e., $g(x_1,...,x_N)\not=\sum_i g_i(x_i)$. 

Funding agencies should therefore tackle the following optimization problem:

\begin{align}
g^\star(X)\equiv\mbox{maximize }&g(x_1,...,x_N),\nonumber\\
\mbox{such that } &\sum_{i=1}^Nx_i=X,\nonumber\\
&X^{-}_j\leq x_j\leq X^{+}_j,j=1,...,N.
\label{optGen}
\end{align}

Even under the assumption that $g(x_1,...,x_N)$ is concave, deterministic and stationary, a blind application of the gradient method will soon lead to trouble. As before, the difficulty stems in estimating the gradient of $g(x_1,...,x_N)$. One way to do so would be to keep the funding of all players but one constant and then compute the difference between the two productivities. For high $N$ this is clearly impractical: even in the absence of statistical fluctuations, proximity to the optimal configuration of funds would only be achieved after $O(N)$ grant calls.

Finally, one could question whether productivity functions exist at all. In the most general case, the productivity of a player at call $k$ could also depend on his/her past success in securing grant funds, i.e., it could be a non-deterministic function, not only of $x_i^k, k$, but also of $x_i^1,...,x_i^{k-1}$. On the other hand, it is possible that the much simpler model (\ref{opt}) already represents an accurate description of the scientific practice. This question cannot be settled by pure mathematical reasoning, but through experimental work, e.g., via pilot research programs.

\subsection{More than one funding agency}
In real life there are several funding bodies at play. Depending on the goals of each funding body, there are different optimization scenarios. If these bodies use the same measure of scientific productivity and their goal is just to increase human knowledge, the best they can do is to create a common budget pool and act as if they were a single funding entity. If they fund completely different areas of research, they can use the policies above independently. If what these bodies fund is pretty much the same, and each of these bodies seeks for recognition, then we enter a complicated game-theoretic problem. One can then divide the funding $x_i$ of scientist $i$ between its sources, i.e., $\sum_{s}x_{i,s}=x_i$, and credit each funding agency $t$ with a proportional amount of the total productivity of each scientist, i.e., $g_{i,t}\equiv\frac{x_{i,t}}{\sum_{s} x_{i,s}}g_i(x_i)$. The goal of each funding agency $t$ would be to maximize $\sum_{i=1}^Ng_{i,t}$, disregarding the performance of all the other agencies.

First of all, as a function of $x_{i,s}$, it is immediate to see that $g_{i,s}$ satisfies $g_{i,s}(0)=0$. One can also prove easily that it is also an increasing function of $x_{i,s}$, since

\be
\frac{\partial g_{i,s}}{\partial x_{i,s}}=\frac{g_i(x_i)}{x_i}-\frac{x_{i,t}}{x_i^2}(g_i(x_i)-g_i'(x_i)x_i)\geq \frac{g_i'(x_i)}{x_i}\geq 0.
\ee
\noindent Here the last inequality follows from eq. (\ref{concavFund}) and $x_i\geq x_{i,s}$.

In addition, assuming $g'''_i(x)\geq 0$ for all $x$ (this is the case, e.g., for the power functions (\ref{powerF})), one can prove that $g_{i,s}$ is also concave. Indeed, note that 

\be
x_i^2\frac{\partial^2g_{i,s}}{\partial x^2_{i,s}}=\frac{x_{i,s}}{x_i}\left(2g_i(x_i)-2g'_i(x_i)x_i+g_i''(x_i)x_i^2\right)-2g_i(x_i)+2g'(x_i)x_i.
\label{deriva}
\ee
Define $h(z)\equiv 2g_i(z)-2g'_i(z)z+g_i''(z)z^2$. It can be verified that $h(0)=0$ and $h'(z)=z^2g_i'''(z)\geq 0$. Hence the term between brackets on the right hand side of eq. (\ref{deriva}) is non-negative for $x_i\geq 0$. Since $\frac{x_{i,s}}{x_i}\leq 1$, we thus have that the right hand side of eq. (\ref{deriva}) is upper bounded by $g_i''(x_i)x_i^2\leq 0$.

Since, as a function of $x_{i,s}$, $\{g_{i,s}\}_{i}$ is concave, the problem of maximizing $\sum_{i=1}^Ng_{i,t}$ for fixed values of $\{x_{i,t}:t\not=s\}_i$ is a convex problem. Moreover, $\{g_{i,s}\}_{i}$ also vanishes at zero, and is increasing, so agency $s$ can apply the grant policy (\ref{new_algo}) to find an orbit of configurations close to the optimal value. The conditions of Theorem \ref{egregius} are also met, and so (under capped funding conditions) agency $s$ can similarly use the grant schemes (\ref{zero_orderGrad}), (\ref{zero_order}) to approximate this maximum. 

All this under the assumption, of course, that the other agencies $t\not=s$ meanwhile keep their budget configurations fixed. If all agencies tried to maximize their total credited productivity at the same time, the system would converge towards a Nash equilibrium, where, of course, $\sum_{s}\sum_{i=1}^Ng_{i,s}$ would not in general coincide with the maximum total productivity achievable. Note that these conclusions also hold when the productivity functions used by the agencies differ, i.e., when the (raw) productivity of scientist $i$ is evaluated differently by each agency.

\section{Conclusion}
\label{conclusion}
In this paper, we have proposed a family of schemes to fund theoretical research. Contrary to the rule in academic funding, these schemes do not rely on a project proposal, but on recent academic performance, as quantified by a given figure of merit. We observed that, once the figure of merit is accepted, the distribution of grant funds becomes an academic problem as opposed to a political issue. 

In this regard, we proposed an algorithm to decide the allocation of funds on each grant call. Under certain idealized assumptions, the algorithm is guaranteed to drive the system, via successive grant calls, to budget distributions maximizing the total scientific productivity. We also introduced alternative schemes, based on the notion of average rates, to tackle scenarios with high statistical fluctuations in the scientific productivity or its evaluation. We explored numerically the performance of the gradient and average rate schemes on real data and compared it with the usual way funding agencies deal with theoretical project proposals.

One of the flaws of the proposed framework for research funding is that, like most others, it may discourage theorists from conducting creative or very original research. Indeed, it is a well-documented fact that creative and unusual ideas usually take time to be accepted by experts \cite{creative}. A creative grant applicant may thus receive a poor evaluation on his/her recent research, thus depriving him/her from a well-deserved funding. A reasonable policy to address this matter, proposed in \cite{trouble}, would be to move researchers with a very high variance in their expert evaluations to an entirely different funding program, perhaps relying on random grant schemes, see \cite{random}.

Most worryingly, our models of scientific productivity are \emph{plagued} with ad hoc assumptions. In order to propose a realistic grant scheme, we need basic information regarding the regular practice of research, information that can only be acquired through experiment. How do productivity functions look like? How are they distributed among theoretical researchers? What is the volatility of expert referee scores? The answers to the questions will teach us whether the research policies presented here work better when applied at the level of individual groups or whole research institutes.

In any case, the purpose of this article is not to provide funding bodies with the ultimate grant scheme, but to contribute to the ongoing \emph{academic} discussion on the problem of research funding. This problem won't be solved by university administrators or politicians. The solution, if it exists, will be reached through the scientific method. Because whenever science comes in, reason and truth follow.

\section*{Acknowledgements}
We acknowledge motivating discussions with C. Brukner, A. Spanu, S. Singh, Z. Wang, N. Villanueva and M. Lewenstein.

%%%%%%%%%%%%%%%%%%%%%%%%%%%%%%%%%%%%%%%%%%%%%%%%%%%%%%%%%

%%%%%%%%%%%%%%%% APPENDIX  

%%%%%%%%%%%%%%%%%%%%%%%%%%%%%%%%%%%%%%%%%%%%%%%%%%%%%%%%%

\begin{appendix}

\section{Computation of $g^\star(X)$ for geometric productivity functions}
\label{solGeom}

Let $g_i(x)$ be of the form (\ref{powerF}). Then the derivative of $g_i(x)$ diverges at $x=0$. This implies that the solution of $\{x_i^\star\}_i$ of problem (\ref{opt}) satisfies $x^\star_i>0$ for all $i$. Under these conditions, $\{x_i^\star\}_i$ can be determined by demanding that any infinitesimal transfer of funds between players $i$ and $j$ $x^\star_i\to \tilde{x}_i=x^\star_i+\delta, x^\star_j\to \tilde{x}_j=x^\star_j-\delta$ should not increase the value of the objective function $\sum_{i=1}^Ng_i(\tilde{x}_i)$. This implies that $g_i(x^\star_i+\delta)-g_i(x^\star_i)+g_j(x^\star_j-\delta)-g_j(x^\star_i)\approx \delta(\frac{dg_i(x^\star_i)}dx_i-\frac{dg_j(x^\star_j)}dx_j)\leq 0$ independently of the sign of $\delta$. This can only be true if, for some $\mu>0$,

\be
\frac{dg_i(x^\star_i)}{dx_i}=\frac{1}{\mu},
\label{gradEqApp}
\ee
\noindent for $i=1,...,N$. See also a detailed discussion of these conditions in Appendix~\ref{app:cosmetic}.

It follows that $x^\star_i=(\alpha_iA_i\mu)^{\frac{1}{1-\alpha_i}}$. The condition $\sum_i x_i=X$ is thus translated to

\be
\sum_i (\alpha_iA_i\mu)^{\frac{1}{1-\alpha_i}}=X.
\label{soluPower1App}
\ee

Given a value of $X$, solving the above equation we can determine the value of $\mu$, whose explicit dependence with $X$ we will express by $\mu(X)$. Once $\mu(X)$ is known, the final total productivity is given by:

\be
g^\star(X)=\sum_i A_i(\alpha_iA_i\mu(X))^{\frac{\alpha_i}{1-\alpha_i}}.
\label{soluPower2App}
\ee

\section{Proof of Theorem \ref{egregius}}
\label{proofApp}

We will first prove the theorem for free funding conditions and functions $g_i(x)$ such that $f_i(x)=\frac{g_i(x)}{x}$ is invertible and ranges in $(0,\infty)$.

Note that for any configuration $\tilde{x}$ satisfying $\sum_i\tilde{x}_i=X$, the quantity $1/\lambda$ in eq. (\ref{limit}) must belong to the interval $[\min_i f_i(\tilde{x}_i), \max_i f_i(\tilde{x}_i)]$. In fact, if $x^\sharp$ is the optimal configuration, i.e., the one satisfying $f_i(x_i^\sharp) =\frac{1}{\lambda}$ for all $i$, for any other configuration $\tilde{x}\neq x^\sharp$, one would have at least two indices $i$ and $j$, such that $\tilde{x}_i > x_i^\sharp$ and $\tilde{x}_j < x_j^\sharp$, since $\sum_i \tilde{x}_i =\sum_i x^\sharp_i=X$. Then $f_i(\tilde{x}_i) < f_i(x_i^\sharp)=\frac{1}{\lambda}$ and $f_j(\tilde{x}_j) > f_j(x_j^\sharp)=\frac{1}{\lambda}$.

On the other hand, $\sum_i g_i(x_i^\sharp)=\sum_i x_i^\sharp f_i(x_i^\sharp)=\frac{\sum_ix^\sharp_i}{\lambda}=\frac{X}{\lambda}$. It follows that

\be
g^\sharp(X)=\sum_i g_i(x^\sharp_i)\in X\left[\min_i f_i(\tilde{x}_i), \max_i f_i(\tilde{x}_i)\right].
\label{bounds}
\ee

Equation (\ref{bounds}) implies that the statement of Theorem \ref{egregius} holds iff, for any feasible distribution of funds $\{x_i\}$, there exists another feasible distribution of funds $\{\bar{x}_i\}_i$ such that

\be
f_i(\bar{x}_i)\geq \frac{1}{2}\frac{\sum_j x_j f_j(x_j)}{X}=:\frac{\bar{f}}{2}.
\label{reforE}
\ee

Let us see why. If Theorem \ref{egregius} is true, then 

\be
g^\sharp(X)=Xf_i(x^\sharp_i)\geq \frac{1}{2}g^\star(X)\geq\frac{1}{2}\sum_j x_j f_j(x_j),
\ee
\noindent for $i=1,...,N$, and for all feasible $\{x_j\}_j$. Dividing by $X$ and identifying $\bar{x}$ with $x^\sharp$, we arrive at Equation (\ref{reforE}).

Conversely, if eq. (\ref{reforE}) holds for $x=x^\star$, then by eq. (\ref{bounds}) we have that

\be
g^\sharp(X)\geq X\min_if_i(\bar{x}_i)\geq \frac{1}{2}\sum_j x^\star_j f_j(x^\star_j)=\frac{g^\star(X)}{2}.
\ee

Assuming that $\{f_i\}_i$ are invertible (and decreasing), we have that equation (\ref{reforE}) is equivalent to:

\be
\bar{x}_i\leq f_i^{-1}\left(\frac{\bar{f}}{2}\right).
\label{condi}
\ee

\noindent Summing on $i$ and taking into account the normalization constraint we arrive at

\be
X\leq \sum_i f_i^{-1}\left(\frac{\bar{f}}{2}\right).
\label{fundamental}
\ee

\noindent Conversely, if the above condition is satisfied, then one can define 

\be
\bar{x}_i=X \frac{f_i^{-1}\left(\frac{\bar{f}}{2}\right)}{\sum_lf_l^{-1}\left(\frac{\bar{f}}{2}\right)}.
\ee
\noindent Then one can verify that $\{\bar{x}_i\}$ satisfy (\ref{condi}) and the normalization constraint. Eq. (\ref{fundamental}) is hence a reformulation of teh statement of Theorem \ref{egregius}.

Let us rewrite eq. (\ref{fundamental}) as

\be
\frac{2}{X\bar{f}}\sum_i x_if_i(x_i)F_i\left(x_i,\frac{\bar{f}}{2}\right)\geq 1,
\label{lessFund}
\ee
\noindent where $F_i(x,y)=\frac{yf_i^{-1}(y)}{g_i(x)}$. Define $p_i\equiv \frac{x_i}{X}$, $\tilde{p}_i\equiv p_i\frac{f_i(x_i)}{\bar{f}}$, and note that both $\{p_i\}_i$ and $\{\tilde{p}_i\}_i$ are normalized probability distributions on the variable $i=1,...,N$.

Now, it can be seen that there exists a non-negative number $f_0$ such that

\be
\sum_{i:f_i(x_i)\geq f_0}\tilde{p}_i, \sum_{i:f_i(x_i)\leq f_0}\tilde{p}_i\geq \frac{1}{2}.
\label{mitad}
\ee
Observe that the second equation implies that

\be
\frac{\bar{f}}{2}\leq \sum_{i:f_i(x_i)\leq f_0}p_if_i(x_i)\leq f_0.
\ee

Putting all together, we have that

\begin{align}
&\frac{1}{X}\sum_i f_i^{-1}\left(\frac{\bar{f}}{2}\right)=2\sum_i\tilde{p}_iF_i\left(x_i,\frac{\bar{f}}{2}\right)\geq\nonumber\\
&2\sum_{i:f_i(x_i)\geq f_0}\tilde{p}_iF_i\left(x_i,\frac{\bar{f}}{2}\right)\geq \nonumber\\
&2\left(\sum_{i:f_i(x_i)\geq f_0}\tilde{p}_i\right)\left(\min_i\left\{F_i\left(x_i,\frac{\bar{f}}{2}\right):f_i(x_i)\geq \frac{\bar{f}}{2}\right\}\right)\geq \nonumber\\
&\min_i\left\{F_i\left(x_i,\frac{\bar{f}}{2}\right):f_i(x_i)\geq \frac{\bar{f}}{2}\right\},
\end{align}
\noindent where the last inequality follows from eq. (\ref{mitad}).

In Appendix \ref{decreas}, we prove that, under the assumption that $g_i(x)$ admits a second derivative, $F_i(x,y)$ is a decreasing function of $y$. This means that, for $f_i(x_i)\geq \frac{\bar{f}}{2}$, $F_i\left(x_i,\frac{\bar{f}}{2}\right)\geq F_i(x_i,f_i(x_i))=1$. This concludes the proof for free funding conditions and productivity functions such that  $f_i(x)$ is invertible and ranges in $(0,\infty)$.

Now, suppose that $f_i(x)$ is not invertible, or doesn't range from $(0,\infty)$, and suppose also that $X^+_i<\infty$. Then, for any $\delta>0$, we can always find a new concave, increasing function $\tilde{g}_i(x)$, with $\tilde{g}_i(0)=0$, and such that

\begin{enumerate}
\item
$\tilde{f}_i(x)$ satisfies the conditions of invertibility and range. 
\item
$\tilde{g}'_i(x)=0$, for $x\geq X^+_i$.
\item
$|g_i(x)-\tilde{g}_i(x)|\leq \delta$ for $x\in[0,X^+_i]$. 
 
\end{enumerate} 
 
\begin{figure}
  \includegraphics[width=12 cm]{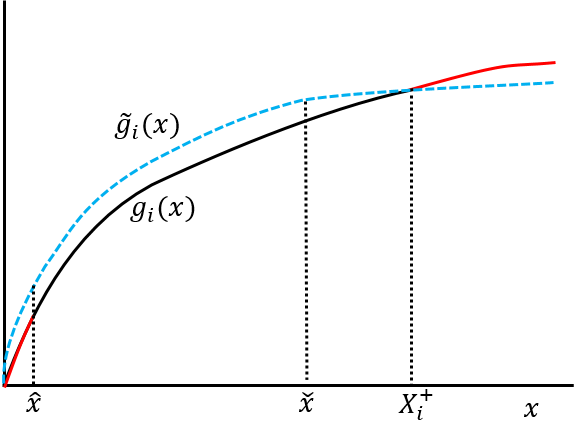}
  \caption{\textbf{Cosmetic surgery.} We modify $g_i(x)+\delta$ from $\hat{x}$ to $0$ such that the new function $\tilde{g}_i(x)$ has an infinite slope at $x=0$. Similarly, we modify $g_i(x)+\delta$ from $\check{x}$ onwards so that the slope of $\tilde{g}_i(x)$ is zero from $X_i^+$ onwards.}
  \label{finalFix}
\end{figure}

Indeed, it suffices to consider the function $\tilde{g}_i(x)=g_i(x)+\delta$, for $x\in[\hat{x},\check{x}]$, with $0<\hat{x}<\check{x}<g_i^{-1}(g_i(X^+_i)-\delta)$. For $x\in [0,\hat{x}]\cup [\check{x},X^+_i]$, once can find a concave, increasing extension $\tilde{g}_i(x)$ of $g_i(x)+\delta$ such that $\tilde{g}_i(x)$ has an infinite slope at $x=0$ and conditions 2,3 above are satisfied. The reader may have a look at Fig~\ref{finalFix} to understand why this is always the case.

Now consider the optimization problem (\ref{opt}) over the productivity functions $\{\tilde{g}_i(x)\}_i$, under free funding conditions. Since the slope of $\tilde{g}_i$ is zero from $x=X^+_i$ onwards, this implies that the optimal solution $\{\tilde{x}^\star_i\}_i$ will satisfy $\tilde{x}_i\star\leq X_i^+$ for all $i$. It therefore coincides with the solution of (\ref{opt}) for capped funding conditions. Since we can choose $\delta>0$ at will, we can do so such that $\sum_i \tilde{g}_i(\tilde{x}^\star_i)\approx g^\star(X)$, where $g^\star(X)$ denotes the optimal solution of the capped problem.

Let $\{\tilde{x}_i\}_i$ be the solution of problem (\ref{optZero}) for the functions $\{\tilde{g}_i\}_i$, assuming free funding conditions. We know, by the previous proof, that 

\be
\sum_i \tilde{g}_i(\tilde{x}_i)\geq\frac{1}{2}\sum_i \tilde{g}_i(\tilde{x}^\star_i)\approx\frac{1}{2}g^\star(X).
\ee

However, in general $\tilde{x}_i\not\leq X^+_i$. Now, define $\tilde{x}'_i=\min(\tilde{x}_i,X^+_i)$. Then it is evident that $\sum_i \tilde{x}_i'\leq\sum_i \tilde{x}_i$ and $\sum_i \tilde{g}_i(\tilde{x}'_i)=\sum_i \tilde{g}_i(\tilde{x}_i)$. The solution $\{\tilde{x}^\sharp_i\}_i$ of the capped problem with the productivity functions $\tilde{g}_i$ will be the result of distributing the excess funds $\sum_i \tilde{x}^\sharp_i-\tilde{x}_i'$ over the players $i$ such that $\tilde{x}_i< X^+_i$. The result can just increase the total productivity, and hence we have that

\be
\sum_i \tilde{g}_i(\tilde{x}^\sharp_i)\gtrsim\frac{1}{2}g^\star(X).
\ee

Finally, it is easy to see that, by decreasing $\delta$, the total productivity of the optimizer of (\ref{optZero}) can be made arbitrarily close to the left hand side of the above equation. It follows that $g^\sharp(X)\geq\frac{1}{2}g^\star(X)$ in the general case.

\section{$F_i(x,y)$ is a decreasing function of $y$}
\label{decreas}

One can easily check that

\be
\frac{dF_i(x,y)}{dy}=\frac{f_i^{-1}(y)+\frac{y}{f_i'(f^{-1}(y))}}{g_i(x)}.
\ee
\noindent Call $z=f^{-1}_i(y)$. Then, written in terms of $g_i(x)$, the numerator of the above equation is proportional to $z\left(1-\frac{g_i(z)}{g_i(z)-xg_i'(z)}\right)$. Now, by eq. (\ref{concavFund}), $g_i(z)-zg_i'(z)$ is non-negative. Since $g_i(z)$ is also non-negative, it hence follows that $\frac{g_i(z)}{g_i(z)-zg_i'(z)}\geq 1$, and so $\frac{dF(x,y)}{dy}\leq 0$, i.e., $F_i(x,y)$ is decreasing on $y$.

\section{Convergence for zero-order method}
\label{app:conv_zero}

\subsection{Conditions on the functions $g_i$}\label{app:cosmetic}
In this section, we will discuss under which conditions the zero-order method will converge to the optimal solution. To simplify the notation, in the following we will assume that the total budget $X$ is normalized, i.e., $X=1$.  This physically corresponds to a change of unit for measuring the budget, so it will not affect the solution. From a mathematical perspective, the same arguments hold for the general case of $X \neq 1$. We will not consider the case of capped funds.

Our assumptions on the single productivity functions $\{g_i\}_i$ are as follows:
\begin{itemize}
\item ${\mathbf{dom}}g_i = [0,1]$,
\item $g_i(0) =0$,
\item $g'_i > 0$, monotonicity (exclude flat case, for uniqueness of solutions)
\item $g''_i < 0$ concavity (exclude linear case).
\end{itemize}

By concavity, it follows that
\begin{equation}\label{eq:conc_te}
g(x_1) + (x_2-x_1)g'(x_1) > g(x_2), \text{ for } x_2>x_1
\end{equation}
which, together with $g(0)=0$ implies
\begin{equation}\label{f>g'}
g(x)- x g'(x) > 0 \Rightarrow \frac{g(x)}{x} > g'(x), \text{ for } x>0.
\end{equation}

The problem in Eq.~\eqref{optZero}, then, becomes
\begin{align}
\mbox{maximize }&\sum_{i=1}^N G_i(x)_i,\nonumber\\
\mbox{such that } &\sum_{i=1}^Nx_i=1,\nonumber\\
& x_j\geq 0,\ j=1,...,N.
\label{optZeroA}
\end{align}
It is convenient to define the functions $f_i(x) := g_i(x)/x$, for all $i$.  We can now derive the conditions for the functions $f_i$ such that the optimal solution $x^\sharp$ for the problem \eqref{optZeroA} satisfies 
\begin{equation}\label{eq:cond_f1}
G'_i(x_i^{\sharp}) =  \frac{g_i(x^\sharp)}{x^{\sharp}} =f_i(x^\sharp)=\frac{1}{\lambda}, \text{ for all } i.
\end{equation}
A necessary condition for a point $x^\sharp$ to be optimal is given by the Karush-Kuhn-Tucker (KKT) conditions \cite{Boyd2004convex}. Moreover, since the problem is concave, with linear inequality constraints and an interior feasible point, by Slater's condition \cite{Boyd2004convex}, the KKT conditions are also sufficient. We can write KKT conditions for the optimal point $x^\sharp$.
\begin{equation}\label{eq:KKT}
\begin{split}
&-\frac{g_i(x_i^\sharp)}{x_i^\sharp}- \mu_i - \nu =0,\\
& \sum_{i=1}^N x_i^\sharp = 1,\\
&x_i^\sharp \geq 0, \text{ for } i=1,\ldots,N\\
&\mu_i \geq 0, \text{ for } i=1,\ldots,N\\
&\mu_i x_i^\sharp = 0, \text{ for } i=1,\ldots,N.
\end{split}
\end{equation}
The last two conditions imply that when $x_i^\sharp > 0$, then $\mu_i = 0$. We want to find conditions on $\frac{g_i(x_i^\sharp)}{x_i^\sharp}= f_i(x_i^\sharp)$, such that there is no solution of Eq.~\eqref{eq:KKT} with $x_i^\sharp=0$. In this case, we can identify $\nu=\frac{1}{\lambda}$ and obtain the condition in Eq.~\eqref{eq:cond_f1}.
 
For example, one could ask that $\lim_{x\rightarrow 0} \frac{g_i(x)}{x}=g'_i(0)=f_i(0)=\infty$, for all $i$. In fact, let us assume that $f_1(0)=\infty$ and $x_1^\sharp = 0$, then to satisfy $\sum_{i=1}^n x_i^\sharp = 1$, at least another $x_i^\sharp$, say $x^\sharp_2$ must be strictly greater than zero. But then $\mu_2 = 0$ and the condition $f_2(x_2^\sharp) = f_1(x_1^\sharp) + \mu_1$ cannot be satisfied for $\mu_1 \geq 0$. More generally, one could simply ask that $f_i$ in zero is ``big enough'' with respect to the other functions $f_j$, $j\neq i$. A sufficient condition to exclude the case $x_i^\sharp = 0$ for some $i$ is given by
\begin{equation}\label{eq:cond_fKKT}
f_i(0) > f_j\left(\frac{1}{N-1}\right), \text{ for all }i, j\neq i.
\end{equation}
This correspond to the configuration in which we assign $0$ to $i$ and an equal amount to $j\neq i$, i.e., $x_i^\sharp = 0$ and $x_j^\sharp=\frac{1}{N-1}$ for $j\neq i$. If $f_i(0)$ is too big, then Eq.~\eqref{eq:KKT} cannot be satisfied. In order to increase the value of some $f_j$ we would have to decrease $x_j^\sharp$, however, given the condition $\sum_{i=1}^N x_i^\sharp = 1$, some other $x^\sharp_{j'}$ should be increased, consequently decreasing the value of $f_{j'}$.

Finally, let us comment on the assumptions on our productivity such as Eq.~\eqref{eq:cond_fKKT}. First, notice that such conditions involve only the local behavior of the function around $x=0$. As a consequence, given any ``actual'' productivity function $g$, we can modify it in a neighborhood of $x=0$ to obtain $\tilde g$ such that, e.g., $\tilde g'(0)=\infty$ but $\tilde g(x)=g(x)$ for all $x > \varepsilon$, for some $\varepsilon > 0$. Applying the iterative method to $g_i$ or $\tilde{g}_i$ for each iterative step $k$ such that $x_i^k > \varepsilon$ will give the same results. Since $\varepsilon$ can be chosen arbitrary small, we can always chose a value such that the values $x_i \leq \varepsilon$ correspond, as a fraction of the total budget $X$, to, e.g., $10^{-3}$ euros. This implies that the difference between $g_i$ and $\tilde{g}_i$ will be relevant only at the step $k$ where we have to redistribute funds of the order of $10^{-3}$ euros. Thus, in practical applications, the assumption on the behavior of $g$ in a neighborhood of $x=0$ implies no loss of generality.

\subsection{Proof of convergence}
We have seen in the previous section that, under condition \eqref{eq:cond_fKKT}, the optimal solution $x^\sharp$ for the problem~\eqref{optZeroA}  satisfies 
\begin{equation}\label{eq:cond_f}
f_i(x^\sharp)=\frac{1}{\lambda}, \text{ for all } i.
\end{equation}
Now, notice that $f'_i(x) = [x g'_i(x) -g_i(x)]/x^2 < 0$, due to Eq.~\eqref{f>g'}, hence  $f$ is monotone and the solution of Eq.~\eqref{eq:cond_f} is unique. In fact, if there were two solutions $x$ and $\tilde{x}$, then $1/\lambda=f_i(x_i)=f_i(\tilde{x}_i) \Rightarrow x_i = f_i^{-1}f_i(\tilde{x}_i)=\tilde{x}_i$ for all $i$.

Moreover, we have seen in Appendix \ref{proofApp} that for any normalized budget distribution $\tilde{x}_i$, we have
\begin{equation}\label{eq:interval}
\frac{1}{\lambda} \in [\min_i f_i(\tilde{x}_i),\max_i f_i(\tilde{x}_i)].
\end{equation}
Given the total productivity $P(x)=\sum_i g_i(x_i)$, the iterative method is defined as 
\begin{equation}\label{eq:iter_app}
x^{k}_i\rightarrow x^{k+1}_i = \frac{g_i(x_i^k)}{P(x^k)},
\end{equation}
with initial point $x^0$ such that $x_i^0 >0$ for all $i$.

To show that the method converges to $x^\sharp$, we will show that for any initial point $x^0$, with $x^0_i > 0$ for all $i$ and $x^0\neq x^\sharp$ the sequences of intervals $I_k := [\min_i f_i(x_i^k), \max_i f_i(x_i^k)]$ satisfies
\begin{equation}
I_{k+1} \subsetneqq I_k, \text{ and } |I_k|\rightarrow 0 \text{ for } k\rightarrow \infty, 
\end{equation}
where for $I=[a,b]$ we define $|I|:= b-a$.

By substituting $g_i(x_i)=x_i f_i(x_i)$ in the definition of Eq.~\eqref{eq:iter_app}, we have
\begin{equation}
\begin{split}
 \frac{f_i(x_i^k)}{P(x^k)} > 1 \Rightarrow {x_i^{k+1}} > {x_i^k};\\
 \frac{f_i(x_i^k)}{P(x^k)} < 1 \Rightarrow {x_i^{k+1}} < {x_i^k};\\
 \frac{f_i(x_i^k)}{P(x^k)} = 1 \Rightarrow {x_i^{k+1}} = {x_i^k} ;
\end{split}
\end{equation}
which implies, by the strict monotonicity of $f_i$ that 
\begin{equation}
\begin{split}
\frac{f_i(x_i^k)}{P(x^k)} > 1 \Rightarrow {f_i(x_i^{k+1})} < {f_i(x_i^k)};\\
 \frac{f_i(x_i^k)}{P(x^k)} < 1 \Rightarrow  {f_i(x_i^{k+1})} > {f_i(x_i^k)};\\
 \frac{f_i(x_i^k)}{P(x^k)} = 1 \Rightarrow  {f_i(x_i^{k+1})} = {f_i(x_i^k)};
\end{split}
\end{equation}
Moreover, by the definition of $P(x^k)$ and $\sum_i x_i = 1$, we have
\begin{equation}
\min_i f_i(x_i^k) \leq P(x^k) = \sum_i x_i^k f_i(x_i^k) \leq \max_i f_i(x_i^k).
\end{equation}

Next we want to prove a condition on the increase (and decrease) for $f_i$ in the iteration, namely
\begin{equation}
\begin{split}
\frac{f_i(x_i^k)}{P(x^k)} > 1 \Rightarrow {f_i(x_i^{k+1})} > {P(x^k)};\\
 \frac{f_i(x_i^k)}{P(x^k)} < 1 \Rightarrow  {f_i(x_i^{k+1})} < {P(x^k)}.
\end{split}
\end{equation}

\noindent Let us consider first the case $0< \alpha < 1$ with $\alpha := f_i(x_i^k)/P(x^k) = x_i^{k+1}/x_i^k$. We have
\begin{equation}
\begin{split}
f_i(x^{k+1}_i) < P(x^k) \Leftrightarrow f_i(\alpha x_i^k) < \frac{f_i(x_i^k)}{\alpha} \Leftrightarrow \alpha \frac{f_i(\alpha x_i^k)}{f_i(x_i^k)} < 1.
\end{split}
\end{equation}
It is sufficient, then, to notice that 
\begin{equation}
\begin{split}
\frac{\alpha f_i(\alpha x)}{f_i(x)}= \alpha \frac{ g_i(\alpha x)}{\alpha x} \frac{x}{g_i(x)} = \frac{g_i(\alpha x)}{g_i(x)}
< 1,
\end{split}
\end{equation}
since $g_i(x)$ is monotonically increasing in $[0,1]$ ($g'_i > 0$), $0<\alpha < 1$ and $x_i^0 > 0$. Analogously, one can prove that in the case $\alpha > 1$,
\begin{equation}
\frac{g_i(\alpha x)}{g_i(x)} > 1 \Rightarrow f_i(x_i^{k+1}) > P(x^k).
\end{equation}

It remains to be proven that $\lim_{k\rightarrow \infty} |I_k| = 0$. We will argue by contradiction. Let us assume that $\lim_{k\rightarrow \infty} I_k = [a,b]$, with $b-a > 0$. Since $\{x^k\}_k$ is bounded, there exist a converging subsequence $\{x^{k_n}\}_n$ with limit $x^{k_n}\rightarrow x^*$, with $\min_i f_i(x_i^*) = a$ and $\max_i f_i(x_i^*) =b$. However, since $b-a>0$, at least one of the following must be true: either $P(x^*) \neq a$ or $P(x^*) \neq b$. Let us assume that $P(x^*) \neq a$, as the other case is identical. Then, by applying the iterative map, we obtain a new interval $I' =[a',b]\subsetneq [a,b]$, in contradiction with the assumption that $[a,b]$ was the limit. $\square$

\subsection{Speed of convergence}
In the following, we will show that the iterative method converges exponentially. First, we need to prove that there exists $\beta < 1$ such that sequence $\{x^k\}_k$ obtained via the iterative method of Eq.~\eqref{eq:iter_app} satisfies
\begin{equation}\label{eq:step_zm}
| f_i(x_i^{k+1}) - P(x^k)| \leq \beta | f_i(x_i^{k}) - P(x^k)|,\ \forall, i,k.
\end{equation}

Let us define $\gamma := P(x^k)^{-1}$. We can rewrite the iterative step as $x_i^{k+1}=g_i(x_i^k)\gamma$ and $f_i(x_i^k)/P(x^k)=\gamma g_i(x_i^k)/ x_i^k$. Let us first assume $f_i(x_i^k)/P(x^k) = 1+ \varepsilon$ with $\varepsilon > 0$, we will treat the other case below. We then have  
\begin{equation}\label{eq:conv1}
\frac{f_i(x_i^{k+1}) - P(x^k)}{f_i(x_i^{k}) - P(x^k)} = \frac{\frac{f_i(x_i^{k+1})}{P(x^k)}-1}{\frac{f_i(x_i^{k})}{P(x^k)}-1}= \frac{\frac{g_i(g_i(x_i^k)\gamma)\gamma}{\gamma g_i(x_i^k)}-1}{\frac{g_i(x_i^k)\gamma}{x_i^k}-1}=\frac{\frac{g_i(g_i(x_i^k)\gamma)}{ g_i(x_i^k)}-1}{\frac{g_i(x_i^k)\gamma}{x_i^k}-1} .
\end{equation}
Let us simplify the expression, using also $\gamma g_i(x_i^k)=(1+\varepsilon) x_i^k$,  the expression \eqref{eq:conv1} becomes
\begin{equation}
\frac{g_i(x_i^k(1+\varepsilon))-g_i(x_i^k)}{\varepsilon g_i(x_i^k)} < \frac{\varepsilon x_i^k g'_i(x_i^k)}{\varepsilon g_i(x_i^k)}= \frac{x_i g'_i(x_i^k)}{g_i(x_i^k)} \leq \beta_{i} < 1,
\end{equation}
where we used Eqs.~\eqref{eq:conc_te},\eqref{f>g'}, respectively, for the two inequalities, and defined $\beta_i = \max_{x\in[0,1]} \frac{x g'_i(x)}{g_i(x)}$. Notice that such a maximum exists since $\frac{x g'_i(x)}{g_i(x)}$ is a continuous function, as it is continuous in $0$, and $[0,1]$ is a closed and bounded interval.

The case $f_i(x_i^k)/P(x^k) = 1- \varepsilon$ with $\varepsilon > 0$ is slightly more complicated. Repeating the initial steps, we obtain
\begin{equation}\label{eq:conv2}
\frac{P(x^k)-f(x_i^{k+1}}{P(x^k)-f(x_i^{k}}=\frac{g_i(x_i^k)-g_i(x_i^k(1-\varepsilon))}{\varepsilon g_i(x_i^k)} < \frac{x_i^k g'_i((1-\varepsilon)x_i^k)}{g(x_i^k)},
\end{equation}
again using Eq.~\eqref{eq:conc_te}. Now, let us drop the indices $i,k$ to make the notation lighter and define
\begin{equation}
H_x(\varepsilon):= \frac{g(x)-g(x(1-\varepsilon))}{\varepsilon g(x)}, \ H_{1x}(\varepsilon) := \frac{x g'((1-\varepsilon)x)}{g(x)}.
\end{equation}
Eq.~\eqref{eq:conv2} becomes $H_x(\varepsilon) < H_{1x}(\varepsilon)$.  We can then verify that the derivative w.r.t. $\varepsilon$ is strictly positive, i.e.,
\begin{equation}
H'_{x}(\varepsilon)=\frac{1}{\varepsilon^2g(x)}\left[-g(x)+g((1-\varepsilon)x)+\varepsilon x g'((1-\varepsilon)x) \right] > 0,
\end{equation}
again using Eq.~\eqref{eq:conc_te}, (compare also to Eq.~\eqref{eq:conv2}). As a consequence, $H_{x}$ is monotonically increasing. Notice that this implies that $H_x$ is continuous in $0$ since it is positive and $\lim_{\varepsilon\rightarrow 0} H_x(\varepsilon)\leq \lim_{\varepsilon\rightarrow 0} H_{1x}(\varepsilon)= xg'(x)/g(x) < 1$. Its maximal value for $\varepsilon \in [0,1]$ is given by $H_x(1)=1$. However, such a value of $\varepsilon$ cannot be reached, since by assumption
\begin{equation}
1-\varepsilon = \frac{f_i(x_i^k)}{P(x^k)} \geq \frac{a_0}{b_0} \Rightarrow \varepsilon \leq 1- \frac{a_0}{b_0},
\end{equation}
where $a_0,b_0$ are the endpoints of the interval $I_0=[a_0,b_0]=[\min_i f_i(x_i^0), \max_i f_i(x_i^0)]$, computed by evaluating all $\{f_i\}_i$ on the first iteration point $x^0$, with $x_i^0 >0$ for all $i$.

We then obtain, for the case $1-\varepsilon$, $\beta_i:=\max_{(x,\varepsilon)\in[0,1]\times[0, 1- a_0/b_0]} H_x(\varepsilon)$. Since, for $x\geq 0$, $H_x(\varepsilon)$ is strictly increasing in $\varepsilon$ and equals $1$ at $\varepsilon=1$, it follows that $\beta_i<1$. Finally, $\beta$ appearing in Eq.~\eqref{eq:step_zm} can be obtained as $\beta:=\max_i \beta_i$.

To complete the proof of exponential speed, we will first show that for each iterative step $k$, and each pair of indices $i,j$ such that $f_i(x_i^k) - P(x^k) > 0$ and  $f_j(x_j^k) - P(x^k) < 0$,
\begin{equation}
|f_i(x_i^{k+1})- f_j(x_j^{k+1})| \leq \beta |f_i(x_i^k) - f_j(x_j^k)|.
\end{equation}
 In fact, $f_i(x_i^k) - P(x^k) > 0 \Rightarrow f_i(x_i^{k+1}) - P(x^k) > 0$ and $f_j(x_j^k) - P(x^k) < 0 \Rightarrow f_j(x_j^{k+1}) - P(x^k) < 0$, hence, we can write
 \begin{align}
 &|f_i(x_i^{k+1})- f_j(x_j^{k+1})| = | f_i(x_i^{k+1}) - P(x^k) | + | P(x^k) - f_j(x_j^{k+1})| 
 \nonumber\\
&\leq \beta ( |f_i(x_i^{k+1}) - P(x^k)| + |P(x^k) - f_j(x_j^{k+1})| )\nonumber\\
&= \beta ( f_i(x_i^{k+1}) - P(x^k) + P(x^k) - f_j(x_j^{k+1}) ) = \beta | f_i(x_i^{k+1}) - f_j(x_j^{k+1}) |.
 \end{align}
Finally, denoting by $m$ the index associated with the minimum at the step $k+1$ i.e., $f_m(x_m^{k+1})=\min_i f_i(x_i^{k+1})$ and $M$ the index associated with the maximum, i.e., $f_M(x_M^{kk+1})=\max_i f_i(x_i^{k+1})$, we can write.
\begin{align}
&|I_{k+1}|=|f_M(x_M^{k+1}) - f_m(x_m^{k+1}) | \leq \beta |f_M(x_M^{k}) - f_m(x_m^{k})| \nonumber\\
&\leq \beta |\max_i f_i(x_i^{k}) - \min_jf_j(x_j^{k})| = \beta |I_{k}|,
\end{align}
which completes the proof of exponential convergence. $\square$

\section{Proof of convergence of the gradient scheme for deterministic productivity functions}
\label{proofConv}
In this Appendix, we will prove that funding policy (\ref{new_algo}), under deterministic, time-independent productivity functions, generates an orbit over the space of budget distributions that stays for most of the time near the optimal productivity. That is, it satisfies eq. (\ref{orbits}).

To do so, we will follow the lines of \cite{subgradient}. We will assume that $|g_i''(x)|\leq \Gamma$ for $x\in [X^-_i,X^+_i]$; that the diameter of the set $B$ of valid budget distributions is $D$; and that $\|\bar{\nabla} g(\bar{x})\|\leq G$, for $\bar{x}\in B$. 

First, by contractivity of projections, we have that

\be
\|\bar{x}^{k+1}-\bar{x}^k\|\leq \|\bar{z}^{k+1}-\bar{x}^k\|=\epsilon\|\tilde{\nabla}^kg\|.
\label{separation}
\ee

Now, $(\tilde{\nabla}^kg)_i=\frac{g_i(x^k_i)-g_i(x^{k-1}_i)}{x^k_i-x^{k-1}_i}=g'_i(x_i^k)+\frac{G^k_i}{2}(x_i^{k-1}-x^k_i)$, where $G^k_i\equiv g''_i(x)$, for some $x\in [x_i^k,x_i^{k+1}]$. It follows that $\tilde{\nabla}^kg(\bar{x}^k)=\bar{\nabla}^k g+\frac{1}{2}\Gamma^k (\bar{x}^{k-1}_i-\bar{x}^{k}_i)$, where $\Gamma^k$ is a diagonal matrix whose (negative) entries are lower bounded by $-\Gamma$. This implies, by eq. (\ref{separation}), that

\be
\|\bar{x}^{k+1}-\bar{x}^k\|\leq \epsilon\left(G+\Gamma\|\bar{x}^k-\bar{x}^{k-1}\|\right)\leq \epsilon \left(G+\frac{1}{2}\Gamma D\right).
\ee

Now, let $\bar{x}^*$ be the budget distribution that maximizes the total scientific productivity. Again, by contractivity of projections, we have that

\begin{align}
&\|\bar{x}^{k+1}-\bar{x}^*\|^2\leq \|\bar{z}^{k+1}-\bar{x}^*\|^2\leq\|\bar{x}^{k}-\bar{x}^*\|^2+\epsilon^2\|\tilde{\nabla}^k g(\bar{x})\|^2+2\epsilon\tilde{\nabla}^kg\cdot(\bar{x}^k-\bar{x}^*)\leq\nonumber\\
&\|\bar{x}^{k}-\bar{x}^*\|^2+\epsilon^2\left(G+\frac{1}{2}\Gamma D\right)^2+2\epsilon\nabla^kg\cdot(\bar{x}^k-\bar{x}^*)+\epsilon\Gamma\|\bar{x}^k-\bar{x}^{k-1}\|\|\bar{x}^k-\bar{x}^\star\|\leq\nonumber\\
&\|\bar{x}^{k}-\bar{x}^*\|^2+\epsilon^2\left(G+\frac{1}{2}\Gamma D\right)\left(G+\frac{3}{2}\Gamma D\right)+2\epsilon\nabla^kg\cdot(\bar{x}^k-\bar{x}^*).
\end{align}
\noindent By induction, we arrive at

\be
\|\bar{x}^{k+1}-\bar{x}^*\|^2\leq \|\bar{x}^{1}-\bar{x}^*\|^2+\epsilon^2k\left(G+\frac{1}{2}\Gamma D\right)\left(G+\frac{3}{2}\Gamma D\right)+2\sum_{j=1}^k\epsilon\bar{\nabla}g\cdot (\bar{x}^k-\bar{x}^*).
\ee
Invoking the inequalities $\|\bar{x}^{k+1}-\bar{x}^*\|^2\geq 0$, $\|\bar{x}^{1}-\bar{x}^*\|^2\leq D^2$ and putting all this together, we have that

\be
\frac{1}{k}\sum_{j=1}^k\bar{\nabla}g\cdot(\bar{x}^*-\bar{x}^j)\leq \frac{1}{2\epsilon k}\left(D^2+\epsilon^2 k \left(G+\frac{1}{2}\Gamma D\right)\left(G+\frac{3}{2}\Gamma D\right) \right).
\label{induction}
\ee

By concavity of $g$, we have that $\bar{\nabla}^kg\cdot(\bar{x}^*-\bar{x}^k)\geq g(\bar{x}^*)-g(\bar{x}^k)\geq 0$. Putting all together, we arrive at

\be
\frac{1}{k}\sum_{j=1}^k g(\bar{x}^*)-g(\bar{x}^j)\leq \frac{1}{2\epsilon k}\left(D^2+\epsilon^2 k \left(G+\frac{1}{2}\Gamma D\right)\left(G+\frac{3}{2}\Gamma D\right)\right).
\ee

In the limit $k\rightarrow \infty$, the right hand side of the equation above can be approximated as $\frac{\epsilon}{2}\left(G+\frac{1}{2}\Gamma D\right)\left(G+\frac{3}{2}\Gamma D\right)$, i.e., it can be made arbitrarily small by decreasing the learning rate $\epsilon$.

\section{Convergence of the average rates scheme A for time-dependent, non-deterministic productivity functions}
\label{convA}

The proof follows the same steps as the convergence of the stochastic subgradient method, see \cite{stoch}. It is also very similar to the proof in Appendix \ref{proofConv}. Call $y^k$ the feasible budget maximizing $\langle G(\bar{x},k)\rangle$, and suppose that $\|y^k-y^{k+1}\|\leq \delta$. Call $(x^j)_j$ the sequence of budgets produced by the average rates scheme A. We will prove that, for high $k$ and suitably chosen learning rate $\epsilon$, with probability $1-O\left(\frac{\delta^{1/4}}{\theta}\right)$, $\frac{1}{k}\sum_{j=1}|g(\bar{x}^j,j)-g(y^j,j)|\leq \theta$.

By Taylor's theorem, we have that $G(x,k)=G(y^k,k)+\bar{\nabla}G(y^k,k)\cdot (x-y^k)+\frac{1}{2}(x-y^k)^T\cdot H\cdot (x-y^k)$, where $H$ is the Hessian of $G(x,k)$ evaluated at a point $c$ within the set $\{p\bar{x}+(1-p)y^k\}$. Such is a diagonal matrix with diagonal elements of value $\frac{g_i'(c_i)-\frac{g_i(c_i)}{c_i}}{c_i}$. By eq. (\ref{concavFund}), we have that each of them is negative. We will assume that, for $i=1,...,N$, there exists a number $h>0$ such that $\frac{1}{2}\left|\frac{g_i'(x)-\frac{g_i(x)}{x}}{x}\right|\geq h$ for all $x\in [X^-_i,X^+_i]$. This can be seen equivalent as taking $g_i(x)$ to be curved at the origin. On the other hand, since $y^k$ is a maximum, we have that $0\geq \frac{dG(tx+(1-t)y^k,k)}{dt}|_{t=0}=\bar{\nabla}G(y^k,k)\cdot (x-y^k)$. This allows us to write 

\be
G(y^k,k)-G(x,k)\geq h\|y^k-x\|^2. 
\label{appOptim}
\ee
\noindent We will use this relation soon. 

Similarly, we will assume that there exists $\gamma>0$ such that $\gamma |g(x,k)-g(y,k)|\leq \|x-y\|$ for all feasible $x,y$. Calling $\bar{\nabla}^k G$ the random vector $(\frac{g_1^k}{x^k_1},\frac{g_2^k}{x^k_2},...)$, we will also assume that $\|\bar{\nabla}^k G\|\leq \Gamma$. Of course, by assumption $\langle\bar{\nabla}^k G\rangle=\bar{\nabla}\langle G(x,k)\rangle$. We will denote by $R$ the radius of the feasible region of budgets.

Now, fix the values of $\{g^j:j=1,...,k-1\}$. Following Appendix \ref{proofConv}, we have that

\be
\|x^{k+1}-y^{k+1}\|\leq \|x^k-y^{k+1}\|^2+2\epsilon\bar{\nabla}^kG\cdot (x^k-y^{k+1})+\epsilon^2\|\bar{\nabla}^kG\|^2.
\ee

In turn, $\|x^k-y^{k+1}\|\leq \|x^k-y^{k}\|+\|y^k-y^{k+1}\|\leq \|x^k-y^{k}\|+\delta$. It follows that $\|x^k-y^{k+1}\|^2\leq \|x^k-y^{k}\|^2+\delta^2+2R\delta$. Also, $\|\bar{\nabla}^kG\|^2\leq \Gamma$. Putting all together, we have that

\begin{align}
&\|x^{k+1}-y^{k+1}\|^2\leq \|x^k-y^{k}\|^2+2\epsilon\bar{\nabla}^kG\cdot (x^k-y^k)+\epsilon^2\Gamma^2+2\epsilon\Gamma\delta+\delta^2+2R\delta=\nonumber\\
&\|x^k-y^{k}\|^2+2\epsilon\bar{\nabla}^kG\cdot (x^k-y^k)+r(\delta,\epsilon),
\end{align}
\noindent with $r(\delta,\epsilon)=2\Gamma \epsilon\delta + \delta^2 + 2R\delta + \Gamma^2\epsilon^2$.

Taking an average over the possible values of $g^{k}$, we have that

\begin{align}
&\langle\|x^{k+1}-y^{k+1}\|^2\rangle_{g^1,...,g^{k-1}}\leq \|x^k-y^{k}\|^2+2\epsilon\bar{\nabla}G\cdot (x^k-y^k)+r(\delta,\epsilon)\nonumber\\
&\leq \|x^k-y^{k}\|^2+2\epsilon (G(x^k,k)-G(y^k,k))+r(\delta,\epsilon).
\end{align}

Now we can fix $\{x^j:j=1,...,k-1\}$ and use the same idea to get rid of the term $\|x^k-y^{k}\|^2$. Iterating, we have that

\be
0\leq \langle\|x^{0}-y^{0}\|^2\rangle\leq R^2+2\epsilon\left\langle\sum_{j=1}^k (G(x^j,j)-G(y^j,j))\right\rangle+kr(\delta,\epsilon).
\ee

Rearranging, we have that

\be
\frac{1}{k}\left\langle\sum_{j=1}^k (G(y^j,j)-G(x^j,j))\right\rangle\leq \frac{r(\delta,\epsilon)}{2\epsilon}+\frac{R^2}{2\epsilon k}.
\ee

Taking the limit $k\to\infty$, we have that the right hand side is bounded by $\frac{r(\delta,\epsilon)}{\epsilon}$. On the other hand, it can be verified that the value of $\epsilon$ that minimizes $\frac{r(\delta,\epsilon)}{\epsilon}$ is $\epsilon^\star=(\sqrt{\delta}\sqrt{\delta + 2 R})/\Gamma$, in which case we have that

\be
\lim_{k\to\infty}\frac{1}{k}\left\langle\sum_{j=1}^k (G(y^j,j)-G(x^j,j))\right\rangle\leq O(\sqrt{\delta}).
\ee

By (\ref{appOptim}), $G(y^j,j)-G(x^j,j)$ can be lower bounded by $h\|y^j-x^j\|^2$, and, in turn, the term $\|y^j-x^j\|$ can be lowerbounded by $\gamma |g(y^j,j)-g(x^j,j)|$. Putting all together, we have that

\begin{align}
&\left(\lim_{k\to\infty}\frac{1}{k}\left\langle\sum_{j=1}^k g(y^j,j)-g(x^j,j)\right\rangle\right)^2\leq \left(\lim_{k\to\infty}\frac{1}{k}\left\langle\sum_{j=1}^k |g(y^j,j)-g(x^j,j)|\right\rangle\right)^2\leq\nonumber\\
&\lim_{k\to\infty}\frac{1}{k}\left\langle\sum_{j=1}^k |g(y^j,j)-g(x^j,j)|^2\right\rangle\leq O(\sqrt{\delta}).
\end{align}

Using the relation $P(Z\geq \theta)\leq \frac{\langle Z\rangle}{\theta}$, valid for any non-negative random variable $Z$, we conclude that 

\begin{align}
&P\left(\left|\lim_{k\to\infty}\frac{1}{k}\sum_{j=1}^k g(y^j,j)-g(x^j,j)\right|\geq \theta\right)\leq\noindent\\
&P\left(\lim_{k\to\infty}\frac{1}{k}\sum_{j=1}^k |g(y^j,j)-g(x^j,j)|\geq \theta\right)\leq O\left(\frac{\delta^{1/4}}{\theta}\right).
\end{align}

\section{Security of the rule of three}
\label{secApp}

In Section \ref{dishonest}, we considered the possibility that Daniel, a member of a large scientific community subject to the rule of three, could win more funds by suitably choosing when to report his research achievements. The purpose of this Appendix is to prove that, in the long run, Daniel cannot expect to obtain more funds than by acting honestly.

Following Table \ref{cheating}, in the $(k-1)^{th}$ call, Daniel's undeclared scientific output equals $g_i^++\sum_{j=0}^{k-1}g_i(x_i^j)-g_i^j$. Multiplying by $\lambda$, invoking the identity $x^{j+1}_i=\lambda g_i^j$ and taking into account that undeclared scientific outputs are non-negative, we have that

\be
\sum_{j=0}^{k} x_i^j\leq x_i^0+\lambda g_i^++\sum_{j=0}^{k-1}\lambda g_i(x_i^j)\leq x_i^0+\lambda g_i^++k\lambda g_i\left(\frac{\sum_{j=0}^{k-1}x_i^j}{k}\right),
\ee
\noindent where the last inequality follows from the concavity of $g_i$.

Define $s^k$ via the relation $\frac{1}{k}\sum_{j=0}^{k-1} x_i^j=x^\sharp_i+s^k$, for $k>0$, and $x_i^0+\lambda g^0=x^\sharp_i+s^0$. Then, the above equation implies

\be
x^\sharp_i+s_{k+1}\leq \frac{1}{k+1}(x^\sharp_i+s^0)+\frac{k}{k+1}\lambda g_i\left(x^\sharp_i+s^{k-1}\right).
\label{step1}
\ee

Now, $\lambda g_i(x^\sharp_i+s)\leq \lambda g_i(x^\sharp_i)+\lambda g'_i(x^\sharp_i)s=x^\sharp_i+\lambda g'_i(x^\sharp_i)s$. In turn, by eq. (\ref{concavFund}), we have that $\lambda g'_i(x^\sharp_i)<\lambda g_i(x^\sharp_i)/x^\sharp_i=1$. Here we have assumed that $x^\sharp_i>0$ and that $g_i$ is curved at the origin. Putting all together, we have that

\be
\lambda g_i(x^\sharp_i+s)\leq x^\sharp_i+\alpha_i s,
\ee
\noindent with $\alpha_i\equiv\lambda g'_i(x^\sharp_i)<1$.

\noindent Applying this relation to the right-hand side of (\ref{step1}) and rearranging, we end up with

\be
s^{k+1}\leq \frac{1}{k+1}s^0+\frac{k}{k+1}\alpha_i s^{k}.
\ee

Since $\alpha_i<1$, it follows, from the above formula, that the sequence $(s^k)_k$ can neither keep growing indefinitely nor converge to a value greater than $0$. This finishes the argument.

\end{appendix}

\bibliography{grantsBib}

\end{document}